%% file: GridImprim.tex
\documentclass[12pt,draft]{amsart}
\usepackage{amsmath,amssymb,bm,enumerate,multicol,multido,pst-3dplot,pstricks,pst-node,caption,subcaption}
\usepackage[mathscr]{euscript}

\newtheorem{theorem}{Theorem}[section]
\newtheorem{proposition}[theorem]{Proposition}
\newtheorem{corollary}[theorem]{Corollary}
\newtheorem{lemma}[theorem]{Lemma}

\theoremstyle{definition}
\newtheorem{definition}[theorem]{Definition}
\newtheorem{example}[theorem]{Example}
\newtheorem{remark}[theorem]{Remark}
\newtheorem{construction}[theorem]{Construction}

\newtheorem{question}[theorem]{Question}

\definecolor{amber}{rgb}{1.0, 0.7, 0.0}
\definecolor{lightblue}{rgb}{0.45, 0.76, 0.98}
\definecolor{pink}{rgb}{0.97, 0.56, 0.65}
\definecolor{malachite}{rgb}{0.04, 0.85, 0.32}
\definecolor{teagreen}{rgb}{0.82, 0.94, 0.75}
\definecolor{darkgreen}{rgb}{0.01, 0.75, 0.24}

\def\D{\mathscr{D}}

\def\PP{\mathscr{P}}
\def\B{\mathscr{B}}

\def\C{\mathscr{C}}
\def\Aut{{\mathrm{Aut}}}
\def\Sym{{\mathrm{Sym}}}

\def\E{\mathscr{E}}
\def\Kbf{\mathbf{K}}
\def\O{\mathscr{O}}
\def\n{\mathbf{n}}
\def\c{\mathrm{c}}

\newcommand{\Des}[2]{\D^{#1}_{#2}}
\newcommand{\DD}[2]{\D(#1,#2)}

\usepackage{anysize}
\marginsize{2.5cm}{2.5cm}{2.5cm}{2.5cm}

\title{Higher-dimensional grid-imprimitive\\ block-transitive designs}

\author[S.H. Alavi]{Seyed Hassan Alavi}
\address{Seyed Hassan Alavi, Department of Mathematics, Faculty of Science, Bu-Ali Sina University, Hamedan, Iran.
}
\email{alavi.s.hassan@basu.ac.ir and  alavi.s.hassan@gmail.com}

\author[C. Amarra]{Carmen Amarra}
\address{Carmen Amarra$^{\rm\lowercase{a}}$, Centre for the Mathematics of Symmetry and Computation, School of Mathematics and Statistics, The University of Western Australia, 35 Stirling Highway, Crawley, 6009 W.A., Australia. {\rm{\emph{P\lowercase{ermanent address}}}}: Institute of Mathematics, University of the Philippines Diliman, C.P. Garcia Avenue, Quezon City 1101, Philippines
}
\email{mcamarra@math.upd.edu.ph}

\author[A. Daneshkhah]{Ashraf Daneshkhah}
\address{Ashraf Daneshkhah, Department of Mathematics, Faculty of Science, Bu-Ali Sina University, Hamedan, Iran.
}
\email{adanesh@basu.ac.ir}

\author{Alice Devillers}
\address{Alice Devillers, Centre for the Mathematics of Symmetry and Computation, School of Mathematics and Statistics, The University of Western Australia, 35 Stirling Highway, Crawley, 6009 W.A., Australia.}
\email{alice.devillers@uwa.edu.au}

\author{Cheryl E. Praeger}
\address{Cheryl E. Praeger, Centre for the Mathematics of Symmetry and Computation, School of Mathematics and Statistics, The University of Western Australia, 35 Stirling Highway, Crawley, 6009 W.A., Australia.}
\email{Cheryl.Praeger@uwa.edu.au}
\date{\today}

\begin{document}

\maketitle

\begin{abstract}
It was shown in 1989 by Delandtsheer and Doyen that, for a $2$-design with $v$ points and block size $k$, a block-transitive group of automorphisms can be point-imprimitive (that is, leave invariant a nontrivial partition of the point set) only if $v$ is small enough relative to $k$. Recently, exploiting a  construction of block-transitive point-imprimitive $2$-designs given by Cameron and the last author, four of the authors studied $2$-designs admitting a block-transitive group that preserves a two-dimensional grid structure on the point set. Here we consider the case where there a block-transitive group preserves a multidimensional grid structure on points. We provide necessary and sufficient conditions for such $2$-designs to exist in terms of the parameters of the grid, and certain `array parameters' which describe a subset of points (which will be a block of the design). Using this criterion, we construct explicit examples of $2$-designs for grids of dimensions three and four, and pose several open questions.

\medskip\noindent
\emph{MathScinet classification (2010):}\quad { 05B05; 05B25; 20B25}

\end{abstract}

\input{pictures-gridimprim}

\section{Introduction}

A seminal result of Delandtsheer and Doyen~\cite{DD} from 1989 showed that, provided the number $v$ of points of a  $2$-$(v,k,\lambda)$ design is large enough relative to the block-size $k$, then a block-transitive group of automorphisms cannot leave invariant a nontrivial partition of the point-set, that is to say, it cannot be point-imprimitive. For smaller values of $v$ relative to $k$, general constructions of block-transitive, point-imprimitive $2$-designs were given by Cameron and the last author \cite{CP93} in 1993. These results inspired numerous efforts over the decades to identify the kinds of point-imprimitive, block-transitive $2$-designs which might exist, in particular seeking to understand the structure of their lattices of invariant point partitions. Several constructions were given for $2$-designs with a block-transitive  subgroup of automorphisms preserving a $2$-dimensional grid structure on the point set \cite{grids22, Linsp09grid,BMV2018,gridsconstr21,CP93,ZZ2,ZZ3,ZZ}.   Here we consider the possibility of a block-transitive group preserving a multi-dimensional grid structure on points.  We show by explicit constructions that this is possible for grids of dimensions three and four, and our work raises several open questions. We say that a permutation group $H$ on a set $\PP$ is \emph{$s$-grid-imprimitive}, where $s\geq2$, if $\PP$ can be identified with an $s$-dimensional grid, that is, a Cartesian product
    \begin{equation} \label{def:P}
    \PP = \E_1 \times \E_2 \times \ldots \times \E_s, \quad \mbox{where the sets $\E_i$ satisfy $|\E_i| = e_i>1$,}
    \end{equation}  
in such a way that $H \leq G$, where 
    \begin{equation}\label{def:G}
        G := \Sym(\E_1) \times \Sym(\E_2) \times \ldots \times \Sym(\E_s) = S_{e_1} \times S_{e_2} \times \ldots \times S_{e_s},
    \end{equation} 
acting naturally in its product action on $\PP$, that is, $(g_1,\dots,g_s) \in G$ maps 
    \begin{equation} \label{act:G}
    (\delta_1,\dots,\delta_s) \to (\delta_1^{g_1},\dots,\delta_s^{g_s}), \ \text{for $ (\delta_1,\dots,\delta_s)\in\PP$}.
    \end{equation}
    
Here is a broad-brush summary of the results in this paper:
    \begin{enumerate}[(i)]
    
    \item We construct an infinite family of block-transitive, $3$-grid-impimitive, $2$-designs (Construction ~\ref{con:s3} and Proposition~\ref{p:Des(3,p)}).
    
    \item We construct  in Example~\ref{ex:s4p2}  a block-transitive $4$-grid-imprimitive $2$-$(v,k,\lambda)$ design with 
    \[
    v = 2^{16} + 2^8 + 1, \ k = 257,\ \text{and}\ \lambda = 7! \cdot 3! \cdot 13! \cdot 241!/(v \cdot 2^{127} \cdot 3^2 \cdot 5 \cdot 7).
    \]
    
    \item We introduce in \eqref{def:chi} the notion of the \emph{array} of a point-subset $B \subseteq \PP$, and in Definition~\ref{def:Des} define an incidence structure $\DD{B}{G}$ on $\PP$ based on $B$. The array is a multi-dimensional version of the $1$-dimensional parameter tuple $\mathbf{x}$ used in \cite[Section 2]{CP93} to construct and analyse point-imprimitive designs based on a single invariant partition.
    
    \item We give in Theorem~\ref{t:2des} necessary and sufficient conditions in terms of the array of $B$ for $\DD{B}{G}$ to be a block-transitive $s$-grid-imprimitive $2$-design. 
    
    \end{enumerate}   

\begin{definition}\label{def:array}
    Let $I := \{1, 2, \ldots, s\}$, and for each $J \subseteq I$, let $\E_J := \prod_{i \in J} \E_i$, with the convention that $\E_\varnothing = \{ \delta_\varnothing \}$ is a singleton set, and let $\pi_J$ be the natural projection from $\PP$ to $\E_J$. For a nonempty subset $B \subseteq \PP$, define the function
    \begin{equation} \label{def:chi}
    \chi_B : \bigcup_{J \subsetneq I} \E_J \rightarrow \mathbb{N}_{\geq 0}, \quad \text{where} \ \forall \ J \subsetneq I \ \text{and} \ \delta_J \in \E_J, \ (\delta_J)\chi_B := \left| \{ \varepsilon \in B \ | \ \varepsilon \pi_J = \delta_J \} \right|.
    \end{equation}
We call $\chi_B$ the \emph{array function} of $B$; and by the \emph{array of $B$} we mean the multi-set of images  $(\delta_J)\chi_B$, for $\delta_J \in \E_J$, $J \subsetneq I$. In particular, $(\varepsilon)\pi_\varnothing = \delta_\varnothing$ for all points $\varepsilon$, and hence $(\delta_\varnothing)\chi_B = |B|$ for all subsets $B$.
\end{definition}

 The incidence structures that are the focus of our study are constructed from a given subset $B\subseteq \PP$ and group $H$ as follows.

\begin{definition} \label{def:Des}
Let $\PP = \prod_{i =1}^s \E_i$ as in \eqref{def:P}, let $B \subseteq \PP$ with $|B| = k$, and let $\binom{\PP}{k}$ denote the set of all $k$-element subsets of $\PP$. Let $G$ be the group defined in \eqref{def:G}. For any $H \leq G$, define the point-block incidence structure $\DD{B}{H}$ by
    \begin{equation} \label{def:D}
    \DD{B}{H} := \big(\PP,\B(B,H)\big), \quad \ \text{where} \ \B(B,H) := B^H = \left\{ B^g \ | \ g \in H \right\}.
    \end{equation}
\end{definition}

The paper \cite{grids22} investigates the case where $s=2$ and $H = G$. In that case, developing ideas from \cite[Section 3]{CP93}, the point set $\PP = \E_1 \times \E_2$ was interpreted as the edge set of the complete bipartite graph $\Kbf_{e_1,e_2}$, and each subset $B$ as a subgraph $B(\Delta)$ of $\Kbf_{e_1,e_2}$. Necessary and sufficient graph-theoretic conditions in order for $\D(B,G)$ to be a $2$-design were obtained. For $s \geq 3$, we could similarly envisage $\PP$ as the set of maximal cliques of the complete multipartite graph $\Kbf_{e_1, \ldots, e_s}$, but it becomes unwieldy to visualise design conditions in terms of graph theoretic constraints. Instead we have chosen to work with the array function of point subsets. In particular, we have identified necessary and sufficient conditions, in terms of the array function $\chi_B$ of $B$, under which $\DD{B}{G}$ is a $2$-design, which are stated in Theorem \ref{t:2des}. Observe that by construction the group $G$ is block-transitive on $\DD{B}{G}$.

\begin{theorem} \label{t:2des}
Let $\PP$ be as in \eqref{def:P}, let $G$ be as in \eqref{def:G}, let $B$ and $\DD{B}{G}$ be as in Definition \emph{\ref{def:Des}}, and let $\chi_B$ be as in \eqref{def:chi}.
    \begin{enumerate}[\rm (a)]
    \item The  incidence structure $\DD{B}{G}$ is a $2$-design if and only if
        \begin{equation} \label{eq:2des}
        \sum_{\delta_J \in \E_J} \left((\delta_J)\chi_B\right)^2 = k + \frac{k(k-1)}{v-1} \left( \bigg( \prod_{i \in J^\c} e_i \bigg) - 1 \right) \quad \text{for all $J \subsetneq I$ with $J\ne \varnothing$},
        \end{equation}
    where $J^\c := I \setminus J$.
    \item Let $H \leq G$ such that $\DD{B}{H}$ is a $2$-$(v,k,\lambda)$ design, where $v=|\PP|=\prod_{i\in I}e_i$ and $k=|B|$. Then 
    \begin{enumerate}[(i)]
        \item $H$ is transitive on $\PP$;
        \item condition \eqref{eq:2des} holds for each proper nontrivial subset $J$ of $I$; and 
        \item ${\lambda = \frac{k(k-1)}{v(v-1)} \cdot |H:H_B|}$, where $H_B$ is the setwise stabiliser of $B$ in $H$.
     \end{enumerate}
    \end{enumerate}
\end{theorem}

 The proof of Theorem~\ref{t:2des} is  derived from  Theorem~\ref{t:n2des} that gives different criteria for a $2$-design which are more difficult to apply in practice (see our comments before Lemma~\ref{lem:sets}).

We show, by describing an explicit construction (see Construction \ref{con:s3}) that there are infinitely many values of $e_1, e_2, e_3$ and $k$ for which conditions \eqref{eq:2des} hold with $s=3$, and hence that there are infinitely many block-transitive, $3$-grid-imprimitive $2$-designs $\DD{B}{G}$. The input to our construction is an integer $p \geq 2$, and the output is a block-transitive, $3$-grid-imprimitive $2$-design. Using a somewhat analogous method, we were able to construct an example of a block-transitive, $4$-grid-imprimitive $2$-design corresponding to the integer $p=2$, and it remains an open problem to extend this construction to arbitrary $p \geq 3$ and/or $s \geq 5$.

\begin{theorem}
There exist infinitely many block-transitive, $3$-grid-imprimitive  $2$-designs, and at least one block-transitive, $4$-grid-imprimitive $2$-design.
\end{theorem}

In \cite[p.~39]{CP93} an example of a block-transitive, $2$-grid-imprimitive $3$-design was given, while in  \cite[Section 5]{grids22} several more examples were constructed  and an approach to searching for and classify them was discussed. However only a finite number of block-transitive, $2$-grid-imprimitive $3$-designs were found, and it was asked  \cite[Problem 1]{grids22} whether or not infinitely many such designs exist. Moreover, we have not found any  $3$-designs which are  block-transitive and $3$-grid-imprimitive.  We summarise our open questions about block-transitive grid-imprimitive designs as follows:

\begin{question}
Do there exist
    \begin{enumerate}[(a)]
    \item infinitely many block-transitive, $4$-grid-imprimitive $2$-designs?
    \item any block-transitive, $s$-grid-imprimitive $2$-designs with $s \geq 5$?
    \item any block-transitive, $3$-grid-imprimitive $3$-designs?
    \end{enumerate}
\end{question}

Using the conditions \eqref{eq:2des}, we are able to obtain necessary conditions in terms of the array function $\chi_B$ of $B$, in order for $\DD{B}{G}$ to be a \emph{flag-transitive $2$-design}. These are presented in Proposition \ref{prop:ft}. None of the designs in Constructions \ref{con:s2} or \ref{con:s3}, or in Example \ref{ex:s4p2}, is flag-transitive, and while infinitely many flag-transitive $2$-grid-imprimitive $2$-designs $\DD{B}{G}$ are known (see \cite[Example 4.4 and Lemma 4.5]{grids22}), the only flag-transitive $2$-design we have found which admits a block-transitive $s$-grid-imprimitive subgroup of automorphisms, with $s\geq3$, is a $2-(16,6,2)$ biplane with $s=3$, see Section~\ref{s:ex}:

\begin{question}
Do there exist any  $2$-designs admitting a flag-transitive, $s$-grid-imprimitive subgroup of automorphisms with $s\geq3$?
\end{question}

\section{Properties of subset arrays and their equivalence} \label{s:array}

In this section, we investigate properties of the array function $\chi_B$ of a point subset $B$. Recall from the definition \eqref{def:chi} that $\chi_B$ is a map with domain $\bigcup_{J \subsetneq I} \E_J$, where $\E_J = \prod_{j \in J} \E_j$. We note the properties mentioned in Definition~\ref{def:array} for $J=\varnothing$.  

The natural product action of $G$ on $\PP$ induces an action of $G$ on $\E_J$, for any $J \subseteq I$, namely,
    \begin{equation} \label{act2:G}
    \text{for any} \ \delta_J \in \E_J \ \text{and} \ g = (g_1, \ldots, g_s) \in G, \quad (\delta_J)^{g} := \big( \delta^{g_j}_j \big)_{j \in J}.
    \end{equation}
It follows, for any $\varepsilon \in \PP$, $J \subseteq I$, and $g \in G$, that $(\varepsilon^g) \pi_J = (\varepsilon \pi_J)^g$. This gives rise  in Lemma \ref{l:chi} to a 
relationship between the array functions $\chi_B$ and $\chi_{B^g}$.

\begin{lemma} \label{l:chi}
Let $\PP$ be as in \eqref{def:P}, let $B \subseteq \PP$, and let $\chi_B$ be as in \eqref{def:chi}. Then, for any $J \subseteq I$, $\delta_J \in \E_J$, and $g \in G$,
    \[
    (\delta_J)\chi_B = ((\delta_J)^g)\chi_{B^g}.
    \]
\end{lemma}

\begin{proof}
By the definition \eqref{def:chi} of the array function $\chi_B$, we have
    \[
    ((\delta_J)^g)\chi_{B^g} = \left|\left\{ \varepsilon \in B^g \ \vline \ \varepsilon \pi_J = (\delta_J)^g \right\}\right|.
    \]
For any $\varepsilon \in B^g$, let $\gamma \in B$ such that $\gamma^g = \varepsilon$. Then $\varepsilon \pi_J = (\gamma^g) \pi_J = (\gamma \pi_J)^g$. So $\varepsilon \pi_J = (\delta_J)^g$ if and only if $(\gamma \pi_J)^g = (\delta_J)^g$, which is in turn equivalent to $\gamma \pi_J = \delta_J$. Thus 
    \[
    ((\delta_J)^g)\chi_{B^g}
    = \left|\left\{ \gamma \in B \ \vline \ \gamma\pi_J = \delta_J \right\}\right|
    = (\delta_J)\chi_B,
    \]
as asserted.
\end{proof}

\begin{definition}\label{d:chi-B-g}
Let $\PP$ be as in \eqref{def:P}. For any $B \subseteq \PP$ and $g \in G$, define the function $\chi^g_B : \bigcup_{J \subsetneq I} \E_J \rightarrow \mathbb{N}_{\geq 0}$ by 
    \[
    \text{for any $J \subsetneq I$ and any $\delta_J \in \E_J$,} \quad (\delta_J)\chi^g_B := \big((\delta_J)^{g^{-1}}\big)\chi_B.
    \]
We shall say that two array functions $\chi_B$ and $\chi_{B'}$ are \emph{$G$-equivalent} if, for some $g \in G$, $\chi_{B'} = \chi_B^g$.
\end{definition}

\begin{corollary} \label{cor:chi}
Let $\PP$ be as in \eqref{def:P}, let $B \subseteq \PP$, and let $\DD{B}{G}$ and $\B(B,G)$ be as in \eqref{def:D}. 
    \begin{enumerate}[\rm (a)]
    \item For $g\in G$, the array function $\chi_{B^g}=\chi_B^g$; and
    \item if $B' \in \B(B,G)$, then $\chi_{B'}$ and $\chi_B$ are $G$-equivalent.
    \end{enumerate}
\end{corollary}

\begin{proof}
(a) Let $g \in G$, $J \subsetneq I$, and $\delta_J \in \E_J$. Then, using Lemma~\ref{l:chi} and Definition~\ref{d:chi-B-g}, 
    \[
    (\delta_J)\chi_{B}^g
    = \big((\delta_J)^{g^{-1}}\big)\chi_{B} 
    = \big((\delta_J)^{g^{-1}}\big)\chi_{(B^{g})^{g^{-1}}}
    = (\delta_J)\chi_{B^g},
    \] 
    and hence $\chi_B^g=\chi_{B^g}$, proving part (a).
    
(b) Let $B' \in \B(B,G)$. Then by \eqref{def:D}, $B' = B^g$ for some $g \in G$, and by part (a), $\chi_{B'} = \chi_{B^g}=\chi_B^g$. Hence $\chi_{B'}$ and $\chi_B$ are $G$-equivalent.
\end{proof}

\begin{remark} \label{rem:B-B'}
The converse of Corollary \ref{cor:chi}(b) is not true in general, that is, there may exist point subsets $B$ and $B'$ which belong to different $G$-orbits but whose array functions are $G$-equivalent, see Example \ref{ex:B-B'}. Recall that the block set of $\DD{B}{G}$ is the $G$-orbit of $B$, so $\DD{B}{G}$ is not necessarily characterised by the arrays of its blocks. Rather, it is possible for two structures $\DD{B}{G}$ and $\DD{B'}{G}$, where $\chi_B$ and $\chi_{B'}$ are $G$-equivalent, to have disjoint block sets. For such $\DD{B}{G}$ and $\DD{B'}{G}$, if $\DD{B}{G}$ is a $2$-design, then it follows from Theorem \ref{t:2des} that $\chi_B$ satisfies the conditions of \eqref{eq:2des}, and hence $\chi_{B'}$ also satisfies these conditions, so $\DD{B'}{G}$ is also a $2$-design. This means moreover that the structure $\big(\PP, \,B^G \cup B'^G\big)$ is also a $2$-design but is not $G$-block-transitive.
\end{remark}

In Example \ref{ex:B-B'}, we illustrate the discussion in Remark~\ref{rem:B-B'}  after introducing the following notation. For each $J \subseteq I$ and $\delta_J \in \E_J$, define $C_{\delta_J}$ to be the set of all points which agree in their $J$-coordinates, that is,
    \begin{equation} \label{def:part}
    C_{\delta_J} := \{ \varepsilon \in \PP \ | \ \varepsilon \pi_J = \delta_J \},
    \end{equation}
and let
    \begin{equation} \label{def:partn}
    \C_J := \big\{ C_{\delta_J} \ \big| \ \delta_J \in \E_J \big\}.
    \end{equation}
Then $\C_J$ is a partition of $\PP$, and in particular $\C_{\varnothing} = \{\PP\}$, and $\C_I = \binom{\PP}{1}$, the partition into singletons. The partition $\C_J$ is nontrivial whenever $\varnothing \neq J \subsetneq I$, and it follows from \eqref{def:chi} and \eqref{def:part} that
    \begin{equation} \label{def2:chi}
    (\delta_J)\chi_B = \big| B \cap C_{\delta_J} \big|.
    \end{equation}
The partition $\C_J$ has $d_J$ classes, each of size $c_J$, where, letting $J^\c := I \setminus J$,
    \begin{equation} \label{def:cd}
    c_J := \prod_{i \in J^\c} e_i \quad \text{and} \quad d_J := \prod_{i \in J} e_i
    \end{equation}
with the convention that the empty products $c_I := 1$ and $d_\varnothing := 1$. Thus $v=|\PP|=c_Jd_J$. Note that if $J' \subsetneq J$, then $\C_J$ is a refinement of $\C_{J'}$, that is, each $\C_{J'}$-class is a disjoint union of $\C_J$-classes. For each $C \in \C_{J'}$, the number of $\C_J$-classes contained in $C$ is
    \[
    \frac{c_{J'}}{c_J} = \frac{\prod_{i \in J'^\c} e_i}{\prod_{i \in J^\c} e_i}
    = \prod_{i \in J \setminus J'} e_i.
    \]

The $G$-action  on $\E_J$ defined in \eqref{act2:G} induces a $G$-action on $\C_J$ by $(C_{\delta_J})^g = C_{(\delta_{J})^g}$, and each partition $\C_J$ is $G$-invariant. Moreover, the group $G$ is the full pointwise stabiliser of this set of partitions, that is to say, $G$ is the largest subgroup of ${\rm Sym}(\PP)$ leaving invariant each of the $\C_J$.

We now give an example to illustrate the comments in Remark \ref{rem:B-B'}.

\begin{example} \label{ex:B-B'}
Let $\PP = \mathbb{Z}_4 \times \mathbb{Z}_4$, $B = \{ 00, 01, 11, 12, 22 \}$, and $B' = \{ 01, 02, 11, 12, 20 \}$. We arrange the points of $\PP$ as a $4 \times 4$ grid $C\times R$, with points in the same column having the same first coordinate, and points in the same row having the same second coordinate. The sets $B$ and $B'$ are illustrated in Figure \ref{fig:B-B'}. Also, as in \cite{grids22}, we interpret the $4\times 4$ grid $C \times R$ as the edge-set of the complete bipartite graph $K_{4,4}$ with vertex bipartition $C\cup R$ such that each of $C$ and $R$ is identified with $\mathbb{Z}_4$, and an edge between $i\in C$ and $j\in R$ is denoted by the pair $ij$. Thus each subset of $C \times R$ corresponds to the edge-set of a subgraph of $K_{4,4}$. In particular, the subgraphs $\Delta$ and $\Delta'$ corresponding to the given subsets $B$ and $B'$  are $\Delta = P_5$ (a path of length $5$) and $\Delta' = C_4 + K_2$ (a disconnected graph with components a cycle $C_4$ of length $4$ and an edge $K_2$), respectively.  Now $B$ and $\Delta$ are as in \cite[Example 4.1]{grids22} with $k = 5$ and $m = 4$, and hence  by \cite[Lemma 4.2]{grids22}, $\DD{B}{G}$ is a $2$-design, where $G = S_4 \times S_4$. Recall that $(\delta_\varnothing)\chi_B = |B| = 4$; and for $J \neq \varnothing$ the array values of $\chi_B$ are listed in Table \ref{tab:B-B'}. Since $I = \{1,2\}$, the only non-empty proper subsets $J$ are $\{1\}$ and $\{2\}$, so we only need to list the values of $(\delta_{\{1\}})\chi_B$ and $(\delta_{\{2\}})\chi_B$. We can deduce these values from the grid on the left in Figure \ref{fig:B-B'} since, by the discussion above, $(\delta_{\{1\}})\chi_B$ is the number of points of $B$ in the column $C_{\delta_{\{1\}}}$ and $(\delta_{\{2\}})\chi_B$ is the number of points of $B$ in the row $C_{\delta_{\{2\}}}$. Similarly, we  obtain the values of $(\delta_J)\chi_{B'}$ from the grid on the right in Figure \ref{fig:B-B'}, and we find that $\chi_B = \chi_{B'}$. So $B$ and $B'$ have $G$-equivalent array values, and it follows from Theorem \ref{t:2des} that $\DD{B'}{G}$ is also a $2$-design. However, since $G\leq \Aut(K_{4,4})$ and $\Delta$ and $\Delta'$ are non-isomorphic subgraphs,  there is no element of $G$ that sends $B$ to $B'$, so $\DD{B}{G}$ and $\DD{B'}{G}$ have disjoint block sets. Finally, as discussed above, the incidence structure $(\PP, \B(B,G) \cup \B(B',G))$ obtained by taking the union of these block sets is also a $2$-design which is not $G$-block-transitive.

   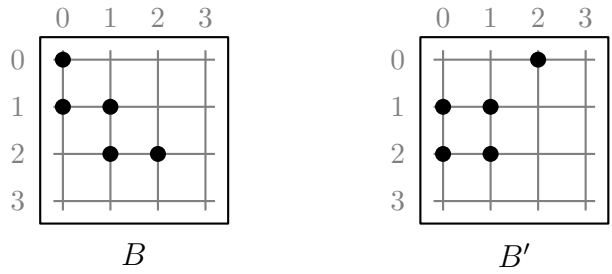
\begin{figure}
        \centering
        \begin{pspicture}(-3,-2)(3,1.5)
        \psset{unit=1.25}
        \exampleblocks
        \end{pspicture}
        \caption{Blocks $B$ and $B'$ in Example \ref{ex:B-B'}}
        \label{fig:B-B'}
    \end{figure}

    \begin{table}[ht]
        \centering
        \begin{tabular}{r|cccc}
        \hline
        $\delta_{\{1\}}$ & $0$ & $1$ & $2$ & $3$ \\
        \hline
        $(\delta_{\{1\}})\chi_B$ & $2$ & $2$ & $1$ & $0$ \\
        \hline
        \end{tabular}
        \quad
        \begin{tabular}{r|cccc}
        \hline
        $\delta_{\{2\}}$ & $0$ & $1$ & $2$ & $3$ \\
        \hline
        $(\delta_{\{2\}})\chi_B$ & $1$ & $2$ & $2$ & $0$ \\
        \hline
        \end{tabular}
        \caption{Array values $(\delta_{\{i\}})\chi_{B}$ for $B$ in Example \ref{ex:B-B'}}
        \label{tab:B-B'}
    \end{table}
\end{example}

\section{Preliminary comments on $\DD{B}{G}$}\label{sec:param}

The incidence structure $\DD{B}{G}$ in Definition \ref{def:Des} is a $1$-design since it admits an automorphism group $G$ which is point-transitive and block-transitive (by construction). The case $s=2$ was studied in \cite{grids22} using a graph theoretic approach. In this paper we consider the case where $s \geq 3$. We first observe that these incidence structures are not $4$-designs; this is a consequence of \cite[Proposition 1.2]{CP93}.

\begin{lemma} \label{l:not4des}
The design $\DD{B}{G}$ in Definition~\emph{\ref{def:Des}} is not a $4$-design.
\end{lemma}

\begin{proof}
Suppose that $\DD{B}{G}$ is a $4$-design. Then, by \cite[Proposition 1.2]{CP93}, the block-transitive group $G$ of automorphisms is $2$-homogeneous on $\PP$, and hence is primitive on $\PP$ by \cite[Lemma 2.30]{PS}. This is a contradiction.
\end{proof}

In order to determine conditions for these designs to be $2$-designs, we will use the following result referred to in \cite{CP93} as a `folklore result'.

\begin{proposition} \cite[Proposition 1.3]{CP93} \label{prop:CP}
Let $t\in\{2,3\}$ and let $H$ be a permutation group on a $v$-element set $\PP$, having orbits $\O_{1}, \ldots, \O_{r}$ on $\binom{\PP}{t}$. Let $B\in \binom{\PP}{k}$, set $\B := \{B^g \mid g \in H\}$. Then for each $i=1,\ldots,r$, there is an integer $\n_i$ such that $\n_i = \big| \binom{B'}{t} \cap \O_i \big|$ for any $B'\in \B$. Moreover, the incidence structure $(\PP, \B)$ is a $t$-design if and only if there exists a constant $c$ such that 
    \[
    \frac{\n_1}{|\O_1|} = \frac{\n_2}{|\O_2|} = \dots = \frac{\n_r}{|\O_r|} = c; \quad \text{moreover,} \ c = \frac{\binom{k}{t}}{\binom{v}{t}} = \frac{\lambda}{|\B|}.
    \]
The group $H$ acts block-transitively on $(\PP,\B)$, and is flag-transitive if and only if the setwise stabiliser $H_B$ of $B$ acts transitively on $B$. 
\end{proposition}

We note that the equality $c = \binom{k}{t}/\binom{v}{t}$ follows from the equations $\n_i = c|\O_i|$, for $i=1,\dots,r$, and the facts that $\sum_i \n_i = \binom{k}{t}$ and $\sum_i |\O_i| = \binom{v}{t}$. The second value for $c$ follows from double-counting the number of incident ($t$-set, block) pairs.

We consider the action of $G$ on $\DD{B}{G}$, and, in particular, we record properties which follow from Proposition~\ref{prop:CP}.

\begin{lemma}\label{lem:flagtr}
Let $\PP$ be as in \eqref{def:P}, let $G$ be as in \eqref{def:G}, let $B$ be a $k$-subset of $\PP$ for some $k \leq |\PP|$, and let $\DD{B}{G}$ be as in Definition~\emph{\ref{def:Des}}. Then $\DD{B}{G}$ is a $1$-$(v,k,\lambda_1)$-design, for some $\lambda_1$, and $G$ is transitive on points and on blocks of $\DD{B}{G}$.
\end{lemma}

\begin{proof}
This follows from the fact that $G$ is point-transitive and block-transitive by Definition~\ref{def:Des}.
\end{proof}

\section{Criteria for $\DD{B}{G}$ to be a $2$-design}

Let $\PP$ and $G$ be as in \eqref{def:P} and \eqref{def:G} and $\DD{B}{G}$ as in Definition~{\ref{def:Des}}, for some subset $B\subseteq \PP$. Then $\DD{B}{G}$ is a $1$-design by Lemma~\ref{lem:flagtr}. We use Proposition~\ref{prop:CP} to derive conditions for $\DD{B}{G}$ to be a $2$-design. These depend on the $G$-orbits on the set $\binom{\PP}{2}$ of unordered pairs of distinct points. There are precisely $2^s-1$ such $G$-orbits, namely, for each non-empty subset $J$ of $I := \{1,2,\ldots,s\}$, the set 
    \begin{equation} \label{def:OJ}
    \O_{J} = \{ \{\delta,\varepsilon\} \mid \text{$\delta$ and $\varepsilon$ differ in exactly their $J$-coordinates} \};
    \end{equation}
 for each such subset $J$ let
    \begin{equation} \label{def:nJ}
    \n_J = \left| \O_{J} \cap {\textstyle\binom{B}{2}}\right|.
    \end{equation}

\begin{theorem}\label{t:n2des}
Let $s \geq 2$, let $\PP, G, B, \DD{B}{G}, I, \O_J, \n_J$ be as above, and let $k=|B|$. Then the following are equivalent:
    \begin{enumerate}[\rm (a)]
    \item $\DD{B}{G}$ is a $2$-design;
    \item for every non-empty subset $J$ of $I$, we have 
        \begin{equation} \label{eq:n2des}
        \n_J = \frac{k(k-1)}{2(v-1)} \prod_{j \in J}(e_j - 1);
        \end{equation}
    \item for all but one non-empty subset $J$ of $I$ condition \eqref{eq:n2des} holds.
    \end{enumerate}
Moreover, $\DD{B}{G}$ is flag-transitive if and only if the setwise stabiliser $G_B$ of $B$ is transitive on $B$.
\end{theorem}

\begin{proof}
First note that, for each $J \subseteq I$, since $G$ is block-transitive, $\n_J = \big|\O_J \cap \binom{B'}{2}\big|$ for each $B' \in \B$, and recall that $|\PP| = v = \prod_{i \in I} e_i$. 

By Proposition~\ref{prop:CP}, $\DD{B}{G}$ is a $2$-design if and only if, for each non-empty subset $J$ of $I$, $\frac{\n_J}{|\O_J|} = \binom{k}{2}/\binom{v}{2} = \frac{k(k-1)}{v(v-1)}$. Moreover, for each $J$, counting ordered pairs $(\delta,\varepsilon)$ with $\{\delta,\varepsilon\}\in\O_J$, we have  $2\cdot|\O_{J}| = v \cdot \prod_{j \in J} (e_j - 1)$.  It follows that $\DD{B}{G}$ is a $2$-design if and only if, for each non-empty subset $J$, $\n_J = \frac{k(k-1)}{2(v-1)} \cdot \prod_{j \in J} (e_j - 1)$. Thus statements (a) and (b) are equivalent. Note that statement (b) involves $2^s - 1$ conditions. Since $\sum_{\varnothing \neq J \subseteq I} \n_J = \binom{k}{2}$, any $2^s - 2$ of these conditions imply the remaining one. Thus statements (a), (b), and (c) are pairwise equivalent.

The condition for flag-transitivity is clear.
\end{proof}

In the next results, we show that a straightforward approach to evaluating the parameters $\n_J$ leads to a double summation involving an alternating 
sum (Lemma~\ref{l:x-nJ}), and therefore a double summation in each of the conditions for verifying we have a $2$-design (Corollary~\ref{c:2des}). Our goal (Theorem~\ref{t:2des}) is to replace these conditions by  more efficiently check-able conditions \eqref{eq:2des} in terms of array values of $B$. To start with, we need the following identities:

\begin{lemma} \label{lem:sets}
For any nonempty finite set $T$, the following identities hold:
    \begin{enumerate}[\rm (a)]
    \item $\sum_{S \subseteq T} (-1)^{|S|} = 0$ \label{binform}
    \item $\prod_{i \in T} (e_i - 1) = \sum_{S \subseteq T} \left( (-1)^{|T \setminus S|} \prod_{i \in S} e_i \right)$ \label{binprod}
    \end{enumerate}
\end{lemma}

\begin{proof}
(a) Let $|T| = n$. Then by the Binomial Theorem~\cite[3.8]{Aig}, $0=(1-1)^n=\sum_{m=0}^n (-1)^m\binom{n}{m}$. Clearly, for each $m=0,\dots,n$, there are precisely $\binom{n}{m}$ subsets $S\subseteq T$ of size $m$, and hence the summand $(-1)^m\binom{n}{m}$ is equal to $\sum_{S \subseteq T, \,|S|=m} (-1)^{|S|}$. Part (a) now follows immediately.

(b) The proof of (b) is similar. Expanding the left hand side $\prod_{i \in T} (e_i - 1) = (e_1 - 1)\dots (e_n - 1)$,  we obtain the sum, over all subsets $S\subseteq T$, of $(-1)^{|T \setminus S|}$ times $\prod_{i \in S} e_i $ (the product of $e_i$ for all $i\in S$ and a factor of $-1$ for each $i\not\in S$).
\end{proof}

For brevity, and when there is no ambiguity, we shall use the notation
    \begin{equation} \label{def:x-chi}
    x_{\delta_J} := (\delta_J)\chi_B \quad \text{for $J \subsetneq I$ and $\delta_J \in \E_J = \prod_{j \in J} \E_j$},
    \end{equation}
where $\chi_B$ is the array function of $B$ defined in \eqref{def:chi}. Hence, in particular, $x_{\delta_\varnothing} = k$, and the array of $B$ is the multiset of all numbers $x_{\delta_J}$, for all $J$ and $\delta_J$, as in \eqref{def:x-chi}. We shall frequently extend this notation  to the case where $J = I$. For this case $\delta_I = \delta \in \PP = \E_I$, and since $\varepsilon \pi_I = \varepsilon$ for any $\varepsilon \in \PP$, $x_{\delta}$ is defined as
    \[
    x_{\delta} := |\{ \varepsilon \in B \ | \ \varepsilon = \delta \}|
    = \begin{cases} 1 &\text{if $\delta \in B$}; \\ 0 &\text{otherwise.} \end{cases}
    \]
Finally, for any $\delta \in \PP$, we define
    \[
    \delta|_J := \delta \pi_J, \quad\text{and note that $\delta|_J\in\E_J$}.
    \]

The following observation will be useful in later proofs. In the sum $\sum_{\delta \in B} x_{\delta|_J}$, the term $x_{\delta|_J}$ occurs as many times as the the number of points in $B$ that belong to the set $C_{\delta|_J}$, and this number is equal to $x_{\delta|_J}$. So $\sum_{\delta \in B} x_{\delta|_J} = \sum_{\bm{\delta} \in B\pi_J} x^2_{\bm{\delta}}$, and since $x_{\bm{\delta}} = 0$ for all $\bm{\delta} \in \E_J \setminus B\pi_J$, we have
    \begin{equation} \label{sum:x-x2}
    \sum_{\delta \in B} x_{\delta|_J}
    = \sum_{\bm{\delta} \in \E_J} x^2_{\bm{\delta}}.
    \end{equation}

\begin{lemma}\label{l:x-nJ}
Let $s \geq 2$ and $I = \{1, 2, \ldots, s\}$, and let $\PP$ be as in \eqref{def:P}. For a given subset $B \subseteq \PP$ and any subset $J \subseteq I$, let $J^\c = I \setminus J$, let the parameter $\n_J$ be as in \eqref{def:nJ}, and for any $\bm{\delta} \in \bigcup_{J \subsetneq I} \E_J$, let $x_{\bm{\delta}}$ be as in \eqref{def:x-chi}. Then
    \begin{equation} \label{eq:x-nJ}
    \n_J 
    = \frac{1}{2} \sum_{S \subseteq J} \sum_{\bm{\delta} \in \E_{J^\c \cup S}} (-1)^{|S|} x_{\bm{\delta}}^2.
    \end{equation}
\end{lemma}

\begin{proof}
By \eqref{def:OJ}, the orbits $\O_J$ are the sets
    \[
    \O_J
    = \{ \{\delta,\varepsilon\} \ | \ \delta_j \neq \varepsilon_j \ \forall\,j \in J \ \text{and} \ \delta_i = \varepsilon_i \ \forall\,i \in J^{\c} \}.    
    \]
It follows from \eqref{def:part} that the condition ``$\delta_i = \varepsilon_i$ for all $i \in J^{\c}$'' implies that $\delta$ and $\varepsilon$ belong to the same $\C_{J^\c}$-part, and the condition ``$\delta_j \neq \varepsilon_j$ for all $j \in J$'' implies that $\delta$ and $\varepsilon$ belong to different $\C_{J^\c \cup \{j\}}$-parts for each $j \in J$. Conversely, if $\delta$ and $\varepsilon$ belong to the same $\C_{J^\c}$-part but different $\C_{J^\c \cup \{j\}}$-parts for $j \in J$, then $\delta_j \neq \varepsilon_j$. Thus
    \begin{equation} \label{OJ-1}
    \O_J
    = \text{$\{\{\delta,\varepsilon\} \ \vline \ \forall\,j \in J, \ \ \delta,\varepsilon \ \text{in same $\C_{J^\c}$-part but in different $\C_{J^\c \cup \{j\}}$-parts}\}$.}
    \end{equation}
Setting  $\O_J^B(\delta) :=  \big\{ \varepsilon \in B \setminus \{\delta\} \ \big| \{ \delta,\varepsilon\} \in \O_J\big\}$, it follows from  \eqref{def:nJ} that 
    \begin{equation}
    \n_J
    = \big| \O_J \cap B^{\{2\}} \big|
    = \frac{1}{2} \cdot \big| \big\{ (\delta,\varepsilon) \ \big| \{ \delta,\varepsilon\} \in \O_J\cap B^{\{2\}} \big\} \big|
    = \frac{1}{2} \sum_{\delta \in B} |\O_J^B(\delta)|. \label{nJ-OJ}
    \end{equation}
Fix $\delta\in B$, and for each $J'\subseteq I$, let $Q_{J'}^B(\delta) := \{ \varepsilon \in B \setminus \{\delta\} \ \big| \delta,\varepsilon \ \text{are in the same $\C_{J'}$-part}\, \}$. Then it follows from \eqref{def:x-chi} that $|Q_{J'}^B(\delta)| = x_{\delta|_{J'}} -1$. Also, from the definition of 
$\O_J^B(\delta)$ and \eqref{OJ-1}, we have 
    \[
    \O_J^B(\delta) = Q_{J^c}^B(\delta)\setminus \left(\bigcup_{j\in J} Q_{J^c\cup\{j\}}^B(\delta)\right). 
    \]
By the Principle of Inclusion and Exclusion \cite[Formula 10.1 in Theorem 10.1]{WvL},
    \[
    \Big| \bigcup_{j \in J}  Q_{J^c\cup\{j\}}^B(\delta) \Big| = 
    \sum_{\varnothing \neq S \subseteq J} (-1)^{|S|+1} \Big| \bigcap_{j \in S} Q_{J^c\cup\{j\}}^B(\delta) \Big|,
    \]
and by the definition of the sets $Q_{J'}^B(\delta)$, we have $\bigcap_{j \in S} Q_{J^c\cup\{j\}}^B(\delta) = Q_{J^c\cup S}^B(\delta)$, which is of cardinality $x_{\delta|_{J^c\cup S}} -1$. Thus 
    \begin{align*}
    \Big| \bigcup_{j \in J}  Q_{J^c\cup\{j\}}^B(\delta) \Big|
    &= \sum_{\varnothing \neq S \subseteq J} (-1)^{|S|+1} \left( x_{\delta|_{J^c\cup S}} -1 \right) \\
    &= -\left( \sum_{\varnothing \neq S \subseteq J} (-1)^{|S|} x_{\delta|_{J^c\cup S}}\right) + \sum_{\varnothing \neq S \subseteq J} (-1)^{|S|} 
    \end{align*}
and by Lemma~\ref{lem:sets}(a), the second summand is equal to $-1$. Since $|Q_{J^c}^B(\delta)|=  x_{\delta|_{J^c}} - 1$, we conclude that
    \[
    \Big| \O_J^B(\delta)\Big| = \left( x_{\delta|_{J^c}} - 1 \right) + \left( \sum_{\varnothing \neq S \subseteq J} (-1)^{|S|} x_{\delta|_{J^c\cup S}}\right) + 1 = \sum_{S \subseteq J} (-1)^{|S|} x_{\delta|_{J^c \cup S}}.
    \]
It follows from this equation and \eqref{nJ-OJ} that
    \begin{equation} \label{eq2:x-nJ}
    \n_J
    = \frac{1}{2} \sum_{\delta \in B} \sum_{S \subseteq J} (-1)^{|S|} x_{\delta|_{J^c \cup S}}
    = \frac{1}{2} \sum_{S \subseteq J} \sum_{\delta \in B} (-1)^{|S|} x_{\delta|_{J^c \cup S}}
    \end{equation}
where, by \eqref{sum:x-x2}, $\sum_{\delta \in B} x_{\delta|_{J^c \cup S}} = \sum_{\bm{\delta} \in \E_{J^\c \cup S}} x_{\bm{\delta}}^2$. The value of $\n_J$ in \eqref{eq:x-nJ} now follows immediately from this equation and \eqref{eq2:x-nJ}.
\end{proof}

\begin{corollary} \label{c:2des}
Let $s \geq 2$ and $I = \{1, 2, \ldots, s\}$, let $\PP$ be as in \eqref{def:P}, let $B$ be a $k$-subset of $\PP$, and let $\DD{B}{G}$ be as in Definition~\emph{\ref{def:Des}}. For any $\bm{\delta} \in \bigcup_{J \subsetneq I} \E_J$ let $x_{\bm{\delta}}$ be as defined in \eqref{def:x-chi}. Then $\DD{B}{G}$ is a $2$-design if and only if 
    \begin{equation} \label{eq2:2des}
    \sum_{S \subseteq J^\c} \sum_{\bm{\delta} \in \E_{J \cup S}} (-1)^{|S|} x_{\bm{\delta}}^2
    = \frac{k(k-1)}{v-1} \prod_{j \in J^\c}(e_j - 1) \quad \text{for all $J \subsetneq I$}.
    \end{equation}
\end{corollary}

\begin{proof}
By replacing $\n_J$ on the left side of \eqref{eq:n2des} with the equivalent expression in \eqref{eq:x-nJ}, and simplifying, we see that the set of conditions in \eqref{eq:n2des} is equivalent to the conditions
    \[
    \sum_{S \subseteq J} \sum_{\bm{\delta} \in \E_{J^\c \cup S}} (-1)^{|S|} x_{\bm{\delta}}^2
    = \frac{k(k-1)}{v-1} \prod_{j \in J}(e_j - 1) \quad \text{for all} \ \varnothing \neq J \subseteq I.
    \]
Note that as $J$ ranges over all nonempty subsets of $I$, its complement $J^\c$ ranges over all proper subsets of $I$. Therefore the above is equivalent to condition \eqref{eq2:2des}.
\end{proof}

Before proving Theorem~\ref{t:2des} we establish a simple fact concerning the equations in \eqref{eq:2des} which is used in the proof, and also in our work in Section~\ref{sec:family}.

\begin{lemma}\label{lem:Jtrivial}
    With notation as in Definition~\textnormal{\ref{def:array}}, the equation in \eqref{eq:2des} holds for $J= \varnothing$ and $J=I$, for any $k$-element subset $B\subseteq \PP = \prod_{i=1}^s\mathbb{Z}_{e_i}$, where $v=|\PP|$. 
\end{lemma}

\begin{proof}
If $J=\varnothing$, then $J^c=I$ so the right hand side of \eqref{eq:2des}  is equal to $k+\frac{k(k-1)}{v-1}(v-1)=k^2$, and the summation on the left hand side of \eqref{eq:2des}  is over the singleton set $\delta_\varnothing$ and  is $((\delta_\varnothing)\chi_B)^2 = |B|^2=k^2$, as noted in the paragraph after \eqref{def:chi}. Thus  \eqref{eq:2des} holds for $J=\varnothing$. Next, if $J=I$, then $J^c=\varnothing$, so the right hand side of \eqref{eq:2des}  is $k+\frac{k(k-1)}{v-1}\cdot 0=k$, and the summation on the left hand side is $\sum_{\delta\in\PP} ((\delta)\chi_B)^2=|B|=k$. Thus  \eqref{eq:2des} also holds for $J=I$.    
\end{proof}

We now prove Theorem \ref{t:2des}, which gives us an alternative set of conditions to \eqref{eq2:2des} that are easier to use in computations.

\begin{proof}[Proof of Theorem \emph{\ref{t:2des}(a)}]
By Corollary \ref{c:2des}, $\DD{B}{G}$ is a $2$-design if and only if condition \eqref{eq2:2des} holds for all $J\subsetneq I$. Thus it is sufficient to show that this set of conditions \eqref{eq2:2des} is equivalent to the set of conditions \eqref{eq:2des} for  $\varnothing\ne J\subsetneq I$. For each $T \subseteq I$, define the notation
    \begin{equation} \label{E_T}
    E_T = \prod_{i \in T} e_i
    \end{equation}
with the convention that $E_T = 1$ when $T = \varnothing$. For an arbitrary subset $J \subseteq I$ it follows from \eqref{def:x-chi} that $\sum_{\delta_J\in\E_J}((\delta_J)\chi_B)^2$ is equal to $\sum_{\delta_J\in\E_J} x^2_{\delta_J}$, so the equation in \eqref{eq:2des} for $J$ is:
    \begin{equation} \label{new6}
    \sum_{\delta_J\in\E_J} x_{\delta_J}^2 = k + \frac{k(k-1)}{v-1} \left( E_{J^\c} - 1 \right).
    \end{equation}
Also for any proper subset $J \subsetneq I$, and any subset $S$ such that $\varnothing\subseteq S \subseteq J^c$, let
    \begin{equation} \label{X_S}
    X_S = (-1)^{|S|} \left( k + \frac{k(k-1)}{v-1} \left( E_{(J \cup S)^\c} - 1 \right) \right).
    \end{equation}
As the sum $\sum_{S \subseteq J^\c} X_S$ will appear repeatedly in the proof below, we begin by proving the following:

\smallskip\noindent
\emph{Claim 1.} For any $J \subsetneq I$,
    \begin{equation} \label{sum-X_S}
    \sum_{S \subseteq J^\c} X_S = \frac{k(k-1)}{v-1} \prod_{i \in J^\c} (e_i - 1).
    \end{equation}

Indeed, we can write
    \begin{align*}
    \sum_{S \subseteq J^\c} X_S
    &= \sum_{S \subseteq J^\c} (-1)^{|S|} \left( k + \frac{k(k-1)}{v-1} \left( E_{(J \cup S)^\c} - 1 \right) \right) \\
    &= k \left(\sum_{S \subseteq J^\c} (-1)^{|S|} \right) + \frac{k(k-1)}{v-1} \left( \sum_{S \subseteq J^\c} (-1)^{|S|} E_{(J \cup S)^\c} - \sum_{S \subseteq J^\c} (-1)^{|S|} \right).
    \end{align*}
Since $J \subsetneq I$, its complement $J^\c$ is nonempty, so that, by Lemma \ref{lem:sets}\eqref{binform}, the sum $\sum_{S \subseteq J^\c} (-1)^{|S|} = 0$. Thus we have 
    \begin{equation} \label{sum2-X_S}
    \sum_{S \subseteq J^\c} X_S
    = \frac{k(k-1)}{v-1} \sum_{S \subseteq J^\c} (-1)^{|S|} E_{(J \cup S)^\c}.
    \end{equation}
For $S\subseteq J^c$, let $S' := (J \cup S)^\c$. Then $S' = J^\c \cap S^\c = J^\c \setminus S$, and as $S$ ranges over all subsets of $J^\c$ so does $S'$. Hence
    \[
    \sum_{S \subseteq J^\c} (-1)^{|S|} E_{(J \cup S)^\c}
    = \sum_{S' \subseteq J^\c} (-1)^{|J^\c \setminus S'|} E_{S'}
    = \prod_{j \in J^\c} (e_j - 1)
    \]
where the last equality follows from Lemma \ref{lem:sets}\eqref{binprod}. Substituting this into \eqref{sum2-X_S} yields \eqref{sum-X_S}. This proves Claim 1.

\smallskip\noindent
{\it Claim 2.}\quad If Condition \eqref{eq:2des} holds for all $J$ such that  $\varnothing\ne J\subsetneq I$, then Condition \eqref{eq2:2des} holds for all $J\subsetneq I$.

Assume that equation \eqref{eq:2des} holds for all subsets $J$ such that $ \varnothing\ne J \subsetneq I$. By Lemma~\ref{lem:Jtrivial}, equation \eqref{eq:2des} also holds for $J = \varnothing$, so these equations hold for all subsets $J \subsetneq I$. It follows that, for any $J$ such that $J \subsetneq I$ and any $S \subseteq J^\c$,
    \[
    \sum_{\bm{\delta} \in \E_{J \cup S}} x^2_{\bm{\delta}}
    = k + \frac{k(k-1)}{v-1} \left( E_{(J \cup S)^\c} - 1 \right).
    \]
Thus, for any $J \subsetneq I$, recalling \eqref{X_S}, we can write
    \begin{equation}
    \sum_{S \subseteq J^\c} \sum_{\bm{\delta} \in \E_{J \cup S}} (-1)^{|S|} x^2_{\bm{\delta}}
    = \sum_{S \subseteq J^\c} (-1)^{|S|} \left( k + \frac{k(k-1)}{v-1} \left( E_{(J \cup S)^\c} - 1 \right) \right)
    = \sum_{S \subseteq J^\c} X_S. \label{eq3:2des}
    \end{equation}
Substituting the expression on the right side of \eqref{sum-X_S} into \eqref{eq3:2des} yields condition \eqref{eq2:2des} for $J$. This proves Claim 2.

\smallskip\noindent
{\it Claim 3.} If Condition \eqref{eq2:2des} holds for  all  $J\subsetneq I$, then Condition \eqref{eq:2des} holds for all $J$ such that  $\varnothing\ne J\subsetneq I$.

Assume that Condition \eqref{eq2:2des} holds for all proper subsets $J\subsetneq I$. We will prove by induction that, for $n=0,1,\dots,s$, the equations in \eqref{eq:2des} hold for all subsets $J$ such that $|J^\c|=n$. This assertion holds for $n=0$ and $n=s$ by Lemma~\ref{lem:Jtrivial}, and so we assume that $1\leq n\leq s-1$ and that the assertion is true for all integers less than $n$.

To complete the proof of the inductive step let $J \subseteq I$ be such that $|J^\c| = n$. Since $1 \leq n \leq s-1$, this means that $\varnothing \ne J \subsetneq I$. For any $S \subseteq J^\c$ let $X_S$ be as defined in \eqref{X_S}. If $S = \varnothing$ then $(J \cup S)^\c = J^\c$, so that $E_{(J \cup S)^\c} = E_{J^\c}$ and
    \begin{align} \label{Xnull}
   X_\varnothing &= k +  \frac{k(k-1)}{v-1} (E_{J^\c} - 1).
    \end{align}
On the other hand, if $\varnothing \neq S \subseteq J^c$, then $J \subsetneq J \cup S$, so the complement $(J \cup S)^\c \subsetneq J^\c$ and $|(J \cup S)^\c| < |J^\c| = n$. Therefore, by the inductive hypothesis, equation \eqref{eq:2des} holds for the set $J \cup S$, that is,
    \[
    \sum_{\bm{\delta} \in \E_{J \cup S}} x_{\bm{\delta}}^2
    = k + \frac{k(k-1)}{v-1} \left( E_{(J \cup S)^\c} - 1 \right)=  (-1)^{|S|} X_S.
    \]
Summing over these subsets we have, 
   \begin{equation}
    \sum_{\varnothing \neq S \subseteq J^\c} \sum_{\bm{\delta} \in \E_{J \cup S}} (-1)^{|S|} x_{\bm{\delta}}^2
    = \sum_{\varnothing \neq S \subseteq J^\c} X_S
    = - X_\varnothing + \sum_{S \subseteq J^\c} X_S. \label{sum3}
    \end{equation}
Thus, substituting into \eqref{sum3} the formulas \eqref{sum-X_S} for $\sum_{S \subseteq J^\c} X_S$ and \eqref{Xnull} for $X_\varnothing$, we get 
    \begin{align}
    \sum_{\varnothing \neq S \subseteq J^\c} \sum_{\bm{\delta} \in \E_{J \cup S}} (-1)^{|S|} x_{\bm{\delta}}^2
    &= - \left( k  + \frac{k(k-1)}{v-1} (E_{J^\c} - 1) \right) + \frac{k(k-1)}{v-1} \prod_{i \in J^\c} (e_i - 1) \notag \\
    &= \frac{k(k-1)}{v-1} \left( \prod_{i \in J^\c} (e_i - 1) - (E_{J^\c} - 1) \right) - k. \label{sum1}
    \end{align}
Note that if $S = \varnothing$ then $\sum_{\bm{\delta} \in \E_{J \cup S}} (-1)^{|S|} x_{\bm{\delta}}^2 = \sum_{\bm{\delta} \in \E_J} x_{\bm{\delta}}^2$. This together with \eqref{sum1} yields
    \begin{align*}
    \sum_{S \subseteq J^\c} \sum_{\bm{\delta} \in \E_{J \cup S}} (-1)^{|S|} x_{\bm{\delta}}^2
    &= \sum_{\bm{\delta} \in \E_J} x_{\bm{\delta}}^2 + \frac{k(k-1)}{v-1} \left( \prod_{i \in J^\c} (e_i - 1) - (E_{J^\c} - 1) \right) -k.
    \end{align*}
Since we assume that \eqref{eq2:2des} holds for $J$, this equality implies, on solving for $\sum_{\bm{\delta} \in \E_J} x_{\bm{\delta}}^2$, that \eqref{new6} holds for $J$, and hence that \eqref{eq:2des} holds for $J$. Thus, by induction, \eqref{eq:2des} holds for all subsets $J \subseteq I$, and Claim 3 is proved.

\medskip
We conclude from Claims 2 and 3 that Condition \eqref{eq:2des} holds for all $J$ such that  $\varnothing\ne J\subsetneq I$ if and only if Condition \eqref{eq2:2des} holds for all $J\subsetneq I$. Thus, by Corollary~\ref{c:2des},  $\DD{B}{G}$ is a $2$-design if and only if Condition \eqref{eq:2des} holds for all $J$ such that  $\varnothing\ne J\subsetneq I$. This proves part (a) of Theorem \ref{t:2des}.
\end{proof}

\begin{proof}[Proof of Theorem \emph{\ref{t:2des}(b)}]
(i) The design $\D(B,H)$ is $H$-block-transitive by construction, so by a result of Block \cite{Block67}, the group $H$ acts transitively on $\PP$. This proves part (i).

(ii) Since $H \leq G$, the block set $B^H \subseteq B^G$. Thus, if $\D(B,H)$ is a $2$-design then $\DD{B}{G}$ is also a $2$-design. It follows from part (a) that condition \eqref{eq:2des} holds for all proper nontrivial subsets $J$ of $I$. Therefore part (ii) holds.

(iii) The values of $v$ and $k$ can be obtained from the definition of $\D(B,H)$. To compute the parameter $\lambda$, recall from Proposition \ref{prop:CP} that $\lambda = \big|B^H\big| \cdot \binom{k}{2}/\binom{v}{2}$. Since $H$ is block-transitive, we have that $\big|B^H\big| = |H:H_B|$, and hence the result follows.
\end{proof}

\section{Examples of $3$-grid imprimitive $2$-designs $\DD{B}{G}$ with small block-size}\label{s:ex}

For $s=3$, Theorem \ref{t:n2des} yields six conditions which together are necessary and sufficient for $\DD{B}{G}$ to be a $2$-design. These six conditions, namely  the ones corresponding to the singleton subsets $\{j\}\subset \{1,2,3\}$, imply the following three divisibility conditions:
\begin{center}
    $2(v-1)$ divides $k(k-1)(e_j-1)$, for $j=1, 2, 3$, where $v=e_1e_2e_3$.
\end{center}
An exhaustive search using {\sc Magma} shows that the  possibilities for $[e_1,e_2,e_3,k]$ with smallest $k$ are $[2,2,4,6]$, $[4,4,4,7]$, $[2,4,7,11]$, and  $[3,3,5,12]$. Again using \textsc{Magma}, we find all examples of $2$-designs $\DD{B}{G}$ with these parameters, up to isomorphism. We list in Table \ref{tab:ex} the values of $\lambda$ and, for each design, a block $B$ from which the block set $\B(B,G)$ can be computed. 

    \begin{table}[ht]
        \centering
        \begin{tabular}{rrrrrrrl}
        \cline{2-8}
        & $e_1$ & $e_2$ & $e_3$ & $v$ & $k$ & $\lambda$ & $B$ \\
        \cline{2-8} \noalign{\vskip\doublerulesep \vskip-\arrayrulewidth} \cline{2-8}
        \small\textsf{1} & $2$ & $2$ & $4$ & $16$ & $6$ & $2$ & $\{ 000, 113, 111, 010, 100, 112 \}$ \\
        \small\textsf{2} &  &  &  &  &  & $6$ & $\{ 000, 002, 103, 001, 012, 112 \}$ \\
        \small\textsf{3} &  &  &  &  &  & $12$ & $\{ 002, 113, 101, 001, 012, 003 \}$ \\
        \small\textsf{4} &  &  &  &  &  & $12$ & $\{ 013, 101, 102, 001, 012, 003 \}$ \\
        \cline{2-8}
        \small\textsf{5} & $4$ & $4$ & $4$ &  & $7$ & $144$ & $\{000, 100, 312, 130, 013, 333, 102\}$ \\
        \small\textsf{6} &  &  &  &  &  &  & $\{000, 100, 012, 122, 311, 301, 121\}$ \\
        \small\textsf{7} &  &  &  &  &  &  & $\{000, 100, 122, 321, 312, 011, 102\}$ \\
        \small\textsf{8} &  &  &  &  &  &  & $\{000, 100, 131, 312, 011, 332, 102\}$ \\
        \cline{2-8}
        \small\textsf{9} & $2$ & $4$ & $7$ &  & $11$ & $4320$ & \parbox[t]{8cm}{\{$000$, $001$, $002$, $003$, $010$, $020$, $034$, $100$, $111$, $125$, $136$\}} \\
        \cline{2-8}
        \small\textsf{10} & $3$ & $3$ & $5$ &  & $12$ & $288$ & \parbox[t]{8cm}{\{$000$, $011$, $012$, $013$, $020$, $101$, $102$, $103$, $111$, $122$, $200$, $213$\}} \\
        \cline{2-8}
        \end{tabular}
        \caption{Examples of $2$-designs $\DD{B}{G}$ with $s=3$ and small $k$}
        \label{tab:ex}
    \end{table}

We make several observations about these data:

\begin{itemize}

\item None of these examples belongs to the infinite family of $2$-designs that we describe in Section \ref{sec:family}, so there is still much to learn about $3$-grid imprimitive block-transitive $2$-designs.

\item Table \ref{tab:ex} lists four pairwise non-isomorphic $2$-designs for which $[e_1,e_2,e_3,k] = [2,2,4,6]$. For these designs $\DD{B}{G}$, the parameters $\n_J$ are $\n_{\{1\}} = \n_{\{2\}} = \n_{\{1,2\}} = 1$ and $\n_{\{3\}} = \n_{\{1,3\}} = \n_{\{2,3\}} = 3$. The example in line 1 has the smallest value for  $\lambda$ and hence  the largest block stabiliser $G_B$ (see Theorem~\ref{t:2des}(iii)); it is a \emph{biplane}, that is, a symmetric $2$-design with $\lambda = 2$. Up to isomorphism there are exactly three $2$-$(16,6,2)$ biplanes,  see \cite{Hus,OR}, and the example in line 1 is the only one which admits a $3$-grid imprimitive, block-transitive subgroup of automorphisms. It is isomorphic to the biplane with point set $\mathbb{Z}^4_2$, and blocks  the set of $\mathbb{Z}^4_2$-translates of $\{ 0000, 1000, 0100, 0010, 0001, 1111 \}$. The full automorphism group $\Aut(\DD{B}{G})$ of this biplane is $2^4.S_6$, while the largest subgroup $H \leq \Aut(\DD{B}{G})$ leaving invariant the $2\times 2\times 4$ grid structure (and permuting the two factors in the $2\times 2$) is permutationally isomorphic to $(S_2\wr S_2)\times S_4$, and contains $G=S_2\times S_2\times S_4$. Although $\Aut(\DD{B}{G})$ acts flag-transitively, a computation with \textsc{Magma} \cite{magma} confirms (unsurprisingly) that the subgroup $H$ is not flag-transitive. We note that this biplane $\DD{B}{G}$ is the only flag-transitive $2$-design we have found which admits a block-transitive, $3$-grid-imprimitive subgroup of automorphisms.

\item Each of the four pairwise non-isomorphic $2$-designs with $[e_1,e_2,e_3,k] = [4,4,4,7]$ has parameters $\n_{\{1\}} = \n_{\{2\}} = \n_{\{3\}} = 1$ and $\n_{\{1,2\}} = \n_{\{1,3\}} = \n_{\{2,3\}} = 3$.  Moreover, for each of them the group $G\leq \Aut(\DD{B}{G})$ is block-regular, that is, $G_B=1$, and  $\lambda = 144$.

\end{itemize}

\section{Necessary conditions for $\DD{B}{G}$ to be flag-transitive}

In this section, we investigate restrictions on the array function $\chi_B$ required for the design  $\DD{B}{G}=(\PP, \B(B,G))$ to be  $G$-flag-transitive. If in fact $\DD{B}{G}$ is $G$-flag-transitive, then, by Theorem \ref{t:n2des}, the block stabiliser $G_B$ is transitive on $B$. It follows that for any $J \subseteq I$, every $\C_J$-class $C$ which intersects $B$ non-trivially contains a constant number of points of $B$, that is to say, using array notation \eqref{def:chi}, 
\begin{center}
    for each $J \subsetneq I$ and each $\delta_J \in \E_J$, $(\delta_J)\chi_{B} = x_{\delta_J} = \big|B \cap C_{(\delta_J)}\big| \in \{0,y_J\}$ 
\end{center}
for some positive integer $y_J$ which depends only on $J$, and we set $y_I := 1$. In this case, we say that the subset $B$ has \emph{uniform array} $(y_J : J \subseteq I)$. Since $\DD{B}{G}$ is $G$-block-transitive, it follows from Corollary \ref{cor:chi} that the array function of all blocks of $\DD{B}{G}$ are $G$-equivalent to  $\chi_{B}$. Thus, if $\DD{B}{G}$ is $G$-flag-transitive, then every block in $\B(B,G)$ also has uniform array $(y_J : J \subseteq I)$. 

First, we examine arithmetical properties of certain synthetically defined quantities $y_J$ which in the $G$-flag-transitive context turn out to be equal to the design parameters $y_J$.

\begin{lemma} \label{lem:arith}
Let $s, k, e_1, \dots, e_s$ be integers, all at least $2$ such that $k < v$, where $v = \prod_{i=1}^s e_i$. Let $D := \gcd(e_1-1, \dots, e_s-1)$, and for each $J \subseteq I$, let $J^c := I \setminus J$, let $c_J$ be as in \eqref{def:cd}, and let 
    \begin{equation} \label{def:y_J}
    y_J := 1 + \frac{k-1}{v-1} \left( c_J - 1 \right) \quad \text{(so in particular $y_\varnothing = k$ and $y_I = 1$).}
    \end{equation}
Assume that $v-1$ divides $(k-1) \cdot D$.
    \begin{enumerate}[\rm (a)]
    \item Then $u := \frac{(k-1)D}{v-1}\in\mathbb{Z}$ and each $y_J$ is a positive integer coprime to $u$. Also for $\varnothing \neq J \subseteq I$ and $j \in J$, the difference $y_{J \setminus \{j\}} - y_J = \frac{k-1}{v-1}(e_j-1) c_J$.
    \item Moreover if, for each $\varnothing \neq J \subseteq I$ and each $j \in J$, the integer $y_J$ divides $\frac{(e_j-1)c_J}{D}$ then, for each $\varnothing \neq J \subseteq I$ and each $j \in J$, $y_J$ divides $y_{J \setminus \{j\}}$ and $1 < \frac{y_{J \setminus \{j\}}}{y_J} < e_j$.
    \end{enumerate}
\end{lemma}

\begin{proof}
By assumption $v-1$ divides $(k-1) \cdot D$, and so  $u := \frac{(k-1)D}{v-1}\in\mathbb{Z}$. Let $\varnothing \neq J \subseteq I$. If $J = I$, then $c_J = c_I = 1$, so that $y_J = y_I = 1$, which is coprime to $u$. On the other hand, if $J \neq I$, then $J^\c =\{ j_1, j_2, \ldots, j_{s-n}\}$, where $|J| = n < s$. Also $c_J = \prod_{i \in J^\c} e_i$, and for any $j \in J^\c$, $c_{J \cup \{j\}}\cdot e_j = \big(\prod_{i \in J^\c \setminus \{j\}} e_i\big)e_j = c_J$. Hence
    \[
    c_J - 1 = c_{J \cup \{j_1\}}(e_{j_1}-1) + \dots + c_{I \setminus \{j_{s-n}\}}(e_{j_{s-n-1}}-1) + (e_{s-n}-1),
    \]
and it follows that $D = \gcd(e_1-1,\dots,e_s-1)$ divides $c_J-1$. Thus we can write \eqref{def:y_J} as
    \[
    y_J = 1 + u \cdot \frac{c_J - 1}{D} \quad \text{for each non-empty subset $J \subseteq I$},
    \]
proving that $y_J\in\mathbb{Z}^+$ and $y_J$ is coprime to $u$. This proves the first assertion of part (a). Now, for $\varnothing \neq J \subseteq I$ and $j \in J$, using the latter equation and noting that $c_{J \setminus \{j\}} = c_Je_j$, we have
    \[ 
    y_{J \setminus \{j\}} - y_J
    = u \cdot \frac{c_{J \setminus \{j\}} - c_J}{D} = u \cdot \frac{(e_j - 1)c_J}{D}
    = \frac{k-1}{v-1}(e_j-1)\prod_{i \in J^\c} e_i. 
    \]
This completes the proof of part (a), noting that the product is `empty' if $J=I$. 

Now assume, in addition, that for each $\varnothing \neq J \subseteq I$ and each $j \in J$ the integer $y_J$ divides $(e_j-1)c_J/D$. Thus, since $u \in \mathbb{Z}$, it follows that $y_J$ divides $y_{J \setminus \{j\}} - y_J$, and therefore $y_J$ divides $y_{J \setminus \{j\}}$, as in part (b). To prove the inequalities on the integer $\frac{y_{J \setminus \{j\}}}{y_J}$, note that $y_{J \setminus \{j\}} > y_J > 1$ by the definition of the $y_J$ in \eqref{def:y_J}, so $\frac{y_{J \setminus \{j\}}}{y_J} > 1$. Also, since $y_J = 1 + u(c_J-1)/D$,
    \[
    y_J e_j
    = e_j + u \cdot \frac{(c_J - 1)e_j}{D}
    = e_j \left( 1 - \frac{u}{D} \right) + u \cdot \frac{c_J e_j}{D}
    = e_j \left( 1 - \frac{u}{D} \right) + u \cdot \frac{c_{J \setminus \{j\}}}{D}.
    \]
Now $k < v$ implies that $\frac{u}{D} = \frac{k-1}{v-1} < 1$, so $1 - \frac{u}{D} > 0$. Then, since $e_j \geq 2$, 
    \[
    e_j \left( 1 - \frac{u}{D} \right) + u \cdot \frac{c_{J \setminus \{j\}}}{D}
    > \left( 1 - \frac{u}{D} \right) + u \cdot \frac{c_{J \setminus \{j\}}}{D} = 1 + u \cdot \frac{c_{J \setminus \{j\}} - 1}{D} = y_{J \setminus \{j\}},
    \]
where the last equality uses the second displayed equation of the proof, so  $\frac{y_{J \setminus \{j\}}}{y_J} < e_j$. 
\end{proof}

We now derive a set of necessary conditions on the uniform array $(y_J : J \subseteq I)$ of a subset $B$ for $\DD{B}{G}$ to be $G$-flag-transitive.

\begin{proposition} \label{prop:ft}
Let $\PP$ be as in \eqref{def:P}, let $G$ be as in \eqref{def:G}, and let $B$ and $\DD{B}{G}$ be as in Definition~\emph{\ref{def:Des}}. Let $D := \gcd(e_1-1, \ldots, e_s-1)$, and for each $J \subseteq I$ let $J^c := I \setminus J$ and let $y_J$ be as in \eqref{def:y_J}. If $\DD{B}{G}$ is a $G$-flag-transitive $2$-design, then the following conditions hold:
    \begin{enumerate}[(i)]
    \item $v-1$ divides $(k-1) \cdot D$, the subset $B$ has  uniform array $(y_J : J \subseteq I)$; and
    \item for each $\varnothing \neq J \subseteq I$, $y_J$ is a positive integer dividing $\frac{D_J}{D} \prod_{i \in J^\c} e_i$, where $D_J := \gcd(e_j - 1 \ | \ j \in J)$.
    \end{enumerate}
Moreover, for each $J \subseteq I$, each class $C \in \C_J$, and each block $B' \in \B$, the intersection size $|B' \cap C| \in \{0,y_J\}$.
\end{proposition}

\begin{proof}
Assume that $\DD{B}{G}$ is a flag-transitive $2$-design. Then by the discussion at the beginning of this section, for each $J \subsetneq I$ and $\delta_J \in \E_J$, we have $(\delta_J)\chi_B \in \{0,y'_J\}$ for some integer $y'_J$ which  is the same for all choices of $\delta_J \in \E_J$, depending only on $J$, and we set $y'_I=1$.

\emph{Claim $1.$\ For each $J \subseteq I$, $y'_J = y_J$ with $y_J$ as in \eqref{def:y_J}.} If $J=I$, then $y'_I=1$, and $y_I=1$ by \eqref{def:y_J}. Also if $J=\varnothing$ then, by Definition~\ref{def:array}, $y'_\varnothing=(\delta_\varnothing)\chi_B=|B|=k=y_\varnothing$ as in \eqref{def:y_J}. Thus we may assume that $\varnothing\ne J\subsetneq I$. Since $\DD{B}{G}$ is a $2$-design, condition \eqref{eq:2des} in Theorem \ref{t:2des} holds for $J$. By \eqref{sum:x-x2} we can write $\sum_{\delta_J \in \E_J} \big((\delta_J)\chi_B\big)^2 = \sum_{\delta \in B} (\delta|_J)\chi_B$, and by the previous paragraph   $(\delta|_J)\chi_B = y'_J$ for all $\delta \in B$ since $\delta|_J\in\E_J$. Thus $\sum_{\delta \in B} (\delta|_J)\chi_B = \sum_{\delta \in B} y'_J = ky'_J$, and by substituting this into the left side of \eqref{eq:2des} and simplifying, we obtain
    \[
    y'_J = 1 + \frac{k-1}{v-1} \left( \Bigg(\prod_{i \in J^\c} e_i\Bigg) - 1 \right) = y_J,
    \]
    where the last equality holds by \eqref{def:y_J}.
This proves Claim 1.

\emph{Claim $2.$ Part (i) holds.} That $B$ has  uniform array $(y_J : J \subseteq I)$ follows from Claim 1. Further, it follows from Claim 1 and Lemma~\ref{lem:arith}(a) that,  for each $J \subseteq I$, $y'_J=y_J$ is a positive integer and thus, by \eqref{def:y_J}, $\frac{k-1}{v-1} \big( \big(\prod_{i \in J^\c} e_i\big) - 1 \big) \in \mathbb{Z}$. Hence $v-1$ divides $(k-1) \cdot \big(\big(\prod_{i \in J^\c} e_i\big) - 1\big)$ for all $J \subseteq I$. As $J$ ranges over all subsets of $I$ so does $J^\c$ and hence, in particular, $v-1$ divides $(k-1) \cdot (e_i - 1)$ for all $i \in I$. Thus $v-1$ divides $(k-1) \cdot D$. Hence part (i) holds, and Claim 2 is proved.

\emph{Claim $3.$ Part (ii) holds.} It follows from Claims 1 and 2 that the hypotheses of Lemma \ref{lem:arith} hold, and in particular $u := \frac{(k-1)D}{v-1} \in \mathbb{Z}^+$. Let $\varnothing \neq J \subseteq I$ and let $j \in J$. Then by Lemma \ref{lem:arith}(a), each $y_J$ is a positive integer coprime to $u$, and
    \[
    y_{J \setminus \{j\}} - y_J
    = \frac{k-1}{v-1} (e_j - 1)c_J
    = u \cdot \frac{e_j - 1}{D} c_J.
    \]
Since each $\C_{J \setminus \{j\}}$-class is a disjoint union of $\C_J$-classes, the number $y_J$ divides $y_{J \setminus \{j\}}$. So $y_J$ divides $y_{J \setminus \{j\}} - y_J = u (e_j - 1) c_J/D$, and since $y_J$ is coprime to $u$, it follows that $y_J$ divides $(e_j - 1)c_J/D$. This is true for all $j \in J$. Hence $y_J$ divides $\gcd\left( (e_j - 1)c_J/D \ | \ j \in J \right) = \gcd( e_j - 1 \ | \ j \in J ) \cdot \frac{1}{D} \cdot c_J = \frac{D_J}{D} \prod_{i \in J^\c} e_i$. Therefore part (ii) holds and Claim 3 is proved.

The final assertion of Proposition~\ref{prop:ft} follows from the fact that $G$-block-transitivity implies that each block has  uniform array $(y_J : J \subseteq I)$.  
\end{proof}

\section{Some infinite families of examples} \label{sec:family}

In this section, we give explicit constructions of infinite families of $2$-designs with point set and block set as in Definition \ref{def:Des}, for $s = 2$ and $s = 3$. In particular, the input for these constructions are two integers $s$ and $p$, both greater than $1$ (with $p$ not necessarily prime), and the output are parameters $e_1, e_2, \ldots, e_s$, a generating block $B$, and a design $\DD{B}{G}$ as in Definition \ref{def:Des} with point set $\PP=\prod_{i=1}^s\mathbb{Z}_{e_i}$. 

\subsection{Recursive approach}\label{s:rec}
In our recursive approach, we construct the $s$-grid-imprimitive $2$-design $\DD{B}{G}$ starting from an $(s-1)$-grid imprimitive $2$-design $\DD{B'}{G'}$, where $G' := S_{e_1} \times \ldots \times S_{e_{s-1}}$, with point set $\PP'=\prod_{i=1}^{s-1}\mathbb{Z}_{e_i}$ and `base block' $B'\subset \PP'$ such that $B$ contains $B'\times\{0\}\subset \PP$, say $B=(B'\times\{0\})\cup B_s$. 

We show in Proposition~\ref{lem:B'} that, provided $B_s$ satisfies certain restrictions, this approach halves the number of conditions which are necessary to be checked to confirm that the structure $\DD{B}{G}$ is a $2$-design.   

We have successfully used this method to construct a new infinite family of $2$-grid-imprimitive $2$-designs $\Des{2}{p}$, for each integer $p>1$ (Construction~\ref{con:s2}), an infinite family of $3$-grid-imprimitive $2$-designs $\Des{3}{p}$ (Construction~\ref{con:s3}) starting from the $2$-grid-imprimitive $2$-design $\Des{2}{p}$ for each integer $p>1$. In addition, we managed to construct a single $4$-grid-imprimitive $2$-design $\Des{4}{2}$ extending the $3$-grid-imprimitive design $\Des{3}{2}$ (Example~\ref{ex:s4p2}). This suggests that it may be possible to extend the example $\Des{4}{2}$ to a family $\Des{4}{p}$ by finding a suitable set $B_4$ that satisfies the conditions of Proposition \ref{lem:B'}\eqref{B_s}.

\begin{question}\label{q:4gi}
    Is there an infinite family of $4$-grid-imprimitive $2$-designs $\Des{4}{p}$ constructible as above from the $3$-grid-imprimitive $2$-designs $\Des{3}{p}$ in Construction~$\ref{con:s3}$?
\end{question}

Throughout this section the parameters $e_1, \ldots, e_s$ are defined as follows:
    \begin{equation} \label{def:e}
    e_1 := p^2 + p + 1, \quad e_i := p^{2^{i-1}} - p^{2^{i-2}} + 1 \ \text{for} \ 2 \leq i \leq s. 
    \end{equation}
Then setting $q_i := p^{2^{i-2}}$, for each $i \geq 2$ we have
    \[ e_i = q_i^2 - q_i + 1 = \frac{q_i^3 + 1}{q_i + 1}. \]
Note that $q_i = q_{i-1}^2$ for each $i \geq 3$. We summarise in Lemma \ref{lem:e} other arithmetic properties of the numbers $e_i$, which will be used frequently in later proofs.

\begin{lemma} \label{lem:e}
Let $p, s > 1$ and let $e_1, \ldots, e_s$ be as in $\eqref{def:e}$.
    \begin{enumerate}[\rm (a)]
    \item \label{e-prod} For any $i \geq 1$,
        \[ \prod_{j=1}^i e_j = \frac{q_{i+1}^3 - 1}{q_{i+1} - 1} = p^{2^i} + p^{2^{i-1}} + 1. \]
    \item \label{e-sum} For any $i \geq 2$,
        \[ \sum_{j=2}^i (e_j - 1) = p^{2^{i-1}} - p \quad \text{and} \quad (p+1) + \sum_{j=2}^i (e_j - 1) = p^{2^{i-1}} + 1. \]
    \item \label{e-quot} For any $i \geq 3$,
        \[ e_i - 1 = (e_{i-1} - 1)\left( q_i + q_{i-1} \right). \]
    \end{enumerate}
\end{lemma}

\begin{proof}
(a) Note that $q_2 = p$, so that $e_1 = q_2^2 + q_2 + 1 = (q_2^3 - 1)/(q_2 - 1)$. 
Thus part (a) holds for $i=1$. If $i \geq 2$ and $\prod_{j=1}^{i-1} e_j = \frac{q_{(i-1)+1}^3 - 1}{q_{(i-1)+1} - 1}$, then inductively, we have
    \[
    \prod_{j=1}^i e_j
    = \left(\prod_{j=1}^{i-1} e_j\right) e_i
    = \frac{q_{i}^3 - 1}{q_{i} - 1} \cdot \frac{q_{i}^3 + 1}{q_{i} + 1}
    = \frac{q_{i}^6 - 1}{q_{i}^2 - 1}
    = \frac{q_{i+1}^3 - 1}{q_{i+1} - 1}
    = q_{i+1}^2 + q_{i+1} + 1
    = p^{2^i} + p^{2^{i-1}} + 1.
    \]
Thus part (a) holds for all $i\geq 1$ by induction.

(b) For the first assertion we have, for each $i\geq2$,
    \[
    \sum_{j=2}^i (e_j - 1)
    = \sum_{j=2}^i (q_j^2 - q_j)
    = \sum_{j=2}^i (q_{j+1} - q_j)
    = q_{i+1} - q_2
    = p^{2^{i-1}} - p.
    \]
The second assertion follows immediately.

(c) For $i\geq3$,
    \begin{align*}
    e_i - 1 
    &= q_i^2 - q_i = q_i^2 - q_{i-1}^2 = (q_i - q_{i-1})(q_i + q_{i-1}) = (q_{i-1}^2 - q_{i-1})(q_i + q_{i-1}) \\
    &= (e_{i-1} - 1)(q_i + q_{i-1}). \qedhere
    \end{align*}
\end{proof}

\subsection{Setting up the inductive step}\label{s:induct}

The incidence structure $\Des{s}{p}$ will have point set $\PP = \prod_{i=1}^s \mathbb{Z}_{e_i}$ and block size $k = p^{2^{s-1}} + 1$ and be defined as an incidence structure $\DD{B}{G}$ as in Definition~{\ref{def:Des}}. To show that $\Des{s}{p}$ is a $2$-design we will need to verify that condition \eqref{eq:2des} holds for each proper nonempty subset $J$ of $I$. Note that $v := \prod_{i=1}^s e_i = p^{2^s} + p^{2^{s-1}} + 1$ by Lemma \ref{lem:e}\eqref{e-prod}, and that
    \[ 
    \frac{k(k-1)}{v-1} = \frac{\left(p^{2^{s-1}} + 1\right)p^{2^{s-1}}}{p^{2^s} + p^{2^{s-1}}} = 1. 
    \]
Thus for each each $\varnothing\ne J\subsetneq I$,  condition \eqref{eq:2des} is equivalent, using \eqref{def:x-chi}, to
    \[ 
    \sum_{\bm{\delta} \in \E_J} x_{\bm{\delta}}^2
    = \left( p^{2^{s-1}} + 1 \right) + \left( \bigg( \prod_{i \in J^\c} e_i \bigg) - 1 \right)
    = p^{2^{s-1}} + \prod_{i \in J^\c} e_i.
    \]
By Theorem \ref{t:2des}, we will need to show that condition \eqref{eq:2des} holds for each $\varnothing\ne J\subsetneq I$. There are $2^s - 2$ such subsets $J$, and hence $2^s - 2$ conditions to check. We show in Proposition  \ref{lem:B'} that, if the generating block $B$ satisfies certain additional hypotheses, then the number of conditions that need to be verified can be reduced by half to $2^{s-1} - 1$.

\begin{proposition} \label{lem:B'}
Let $s, p > 1$ and let the numbers $e_i$, $1 \leq i \leq s$, be as in $\eqref{def:e}$. Let $\PP = \prod_{i=1}^s \mathbb{Z}_{e_i}$ and let $B \subseteq \PP$ be the disjoint union $B = (B' \times \{0\}) \cup B_s$, where $B'\subset \PP' = \prod_{i=1}^{s-1} \mathbb{Z}_{e_i}$ and $B_s\subset \PP$ have the following properties:
    \begin{enumerate}[\rm (a)]
    \item \label{k'}  $|B'| = p^{2^{s-2}} + 1$, and the incidence structure $\DD{B'}{G'}$ with point set $\PP' = \prod_{i=1}^{s-1} \mathbb{Z}_{e_i}$ and $G'=\prod_{i=1}^{s-1}S_{e_i}$  is a $2$-design;
    \item \label{B_s} $B_s$ comprises exactly one element of the $\C_{\{s\}}$-part $C_{\bm{\delta}}$, for each $\bm{\delta} \in \mathbb{Z}_{e_s}\setminus\{0\}$; in particular $|B_s|=e_s-1$.
    \end{enumerate}
Then $|B|=p^{2^{s-1}}+1$, and $\DD{B}{G}$ is a $2$-design if and only if conditions \eqref{eq:2des} hold for all subsets $J \subsetneq I$ such that $s \notin J$.
\end{proposition}

\begin{proof}
Assume that $\PP$ and $B$ are as described and that hypotheses \eqref{k'}--\eqref{B_s} hold. Then  $|B_s|=e_s-1$, so by \eqref{def:e}, $k:=|B|=|B'|+|B_s|= \big(p^{2^{s-2}} + 1\big)+(e_s-1) = p^{2^{s-1}}+1$. If $\DD{B}{G}$ is a $2$-design, then it follows from Theorem \ref{t:2des}(a) that condition \eqref{eq:2des} holds  for all proper nonempty subsets $J \subseteq I$ such that $s \notin J$, and \eqref{eq:2des} also holds for $J=\varnothing$ by Lemma~\ref{lem:Jtrivial}. 

Suppose conversely that condition \eqref{eq:2des} holds for all proper subsets $J\subset I$ that do not contain $s$. To prove that $\DD{B}{G}$ is a $2$-design, we show that condition \eqref{eq:2des} also holds for all proper subsets of $I$ containing $s$. Let $J$ be such a subset, and let $J' := J \setminus \{s\}$. 
Then $\E_J = \E_{J'} \times \E_s$ is the disjoint union $X_\# \cup X_0$, where $X_\# := \E_{J'} \times (\E_s \setminus \{0\}) = \{ \bm{\delta} \in \E_J \ | \ \bm{\delta}\pi_{\{s\}} \neq 0 \}$ and $X_0 := \E_{J'} \times \{0\} = \{ \bm{\delta} \in \E_J \ | \ \bm{\delta}\pi_{\{s\}} = 0 \}$. Thus we can write
    \begin{equation} \label{sum:X+X}
    \sum_{\bm{\delta} \in \E_J} x_{\bm{\delta}}^2
    = \sum_{\bm{\delta} \in X_\#} x_{\bm{\delta}}^2 + \sum_{\bm{\delta} \in X_0} x_{\bm{\delta}}^2,
    \end{equation}
and we evaluate each summand separately. 

\smallskip\noindent
\emph{Claim 1:} The first summand of \eqref{sum:X+X} satisfies  $\sum_{\bm{\delta} \in X_\#} x_{\bm{\delta}}^2
    = p^{2^{s-1}} - p^{2^{s-2}}.$

\smallskip
It follows from the definition of $X_\#$ that (with $C_{\bm{\delta}}$ and $C_{\bm{\gamma}}$ as in \eqref{def:part} for $\bm{\delta}\in\E_J$ and $\bm{\gamma}\in\E_s$, respectively)
    \[
    \bigcup_{\bm{\delta} \in X_\#} C_{\bm{\delta}}
    = \PP' \times (\E_s \setminus \{0\})
    = \{ \varepsilon \in \PP \ | \ \varepsilon_s \neq 0 \}
    = \bigcup_{\bm{\gamma} \in \E_s \setminus \{0\}} C_{\bm{\gamma}}. \]
Since $s \in J$, the partition $\C_J$ in \eqref{def:partn} is a refinement of the partition $\C_{\{s\}}$, so each $\C_J$-part $C_{\bm{\delta}}$ with $\bm{\delta} \in X_\#$ is contained in the $\C_{\{s\}}$-part $C_{\bm{\gamma}}$ with $\bm{\gamma}=\bm{\delta}\pi_{\{s\}} \in \E_s \setminus \{0\}$. (Note that $\C_{J}$ and $\C_{\{s\}}$ are equal in the case $J=\{s\}$.) Therefore $x_{\bm{\delta}} = |B \cap C_{\bm{\delta}}| \leq |B \cap C_{\bm{\gamma}}| = x_{\bm{\gamma}}$. By the definition of $B$, the subset  $B_s$ is $B \cap (\bigcup_{\bm{\delta} \in X_\#} C_{\bm{\delta}})$, so by condition \eqref{B_s}, for any $\bm{\gamma} \in \E_s \setminus \{0\}$, the parameter $x_{\bm{\gamma}} = |B \cap C_{\bm{\gamma}}| = |B_s \cap C_{\bm{\gamma}}| = 1$. Therefore $x_{\bm{\delta}} \leq 1$, so in particular $x_{\bm{\delta}} \in \{0,1\}$ and $x_{\bm{\delta}}^2 = x_{\bm{\delta}}$. Therefore since the sets $B \cap C_{\bm{\delta}}$ are pairwise disjoint, we have $\bigcup_{\bm{\delta} \in X_\#}(B \cap C_{\bm{\delta}}) = B \cap \big( \bigcup_{\bm{\delta} \in X_\#} C_{\bm{\delta}} \big) = B \cap \big( \bigcup_{\bm{\gamma} \in \E_s \setminus \{0\}} C_{\bm{\gamma}} \big) = B_s$ and
    \[
    \sum_{\bm{\delta} \in X_\#} x_{\bm{\delta}}^2
    = \sum_{\bm{\delta} \in X_\#} x_{\bm{\delta}}
    = \sum_{\bm{\delta} \in X_\#} |B \cap C_{\bm{\delta}}|
    = \bigg|\bigcup_{\bm{\delta} \in X_\#}(B \cap C_{\bm{\delta}})\bigg|
    = |B_s|
    = e_s - 1.
    \]
Recalling that $e_s = p^{2^{s-1}} - p^{2^{s-2}} + 1$ by \eqref{def:e}, we conclude that $ \sum_{\bm{\delta} \in X_\#} x_{\bm{\delta}}^2
    = p^{2^{s-1}} - p^{2^{s-2}}$, proving Claim 1.
    
\smallskip\noindent
\emph{Claim 2:} The second summand of \eqref{sum:X+X} satisfies  $\sum_{\bm{\delta} \in X_0} x_{\bm{\delta}}^2
    = p^{2^{s-2}} + \prod_{i \in J^\c} e_i.$

\smallskip
Let $I' := I \setminus \{s\}$. It follows from the definition of $X_0$ and \eqref{def:partn} that
    \[
    \bigcup_{\bm{\delta} \in X_0} C_{\bm{\delta}}
    = \PP' \times \{0\}
    = \{ \varepsilon \in \PP \ | \ \varepsilon_s = 0 \}
    =: C_{0_s} \in \C_{\{s\}}.
    \]
Therefore the restriction to $C_{0_s}$ of the projection map $\pi_{I'} : \PP \rightarrow \PP'$ is a bijection $C_{0_s} \rightarrow \PP'$. Since the sets $C_{\bm{\delta}}$ (for $\bm{\delta}\in X_0$) are pairwise disjoint, the set $\mathscr{F} := \{ C_{\bm{\delta}} \ | \ \bm{\delta} \in X_0 \}$ forms a partition of the $\C_{\{s\}}$-part $C_{0_s}$. We see from the definition of $B$ that $B \cap C_{0_s} = B' \times \{0\}$, so the bijection $\pi_{I'}|_{C_{0_s}}$ sends $B \cap C_{0_s} \rightarrow B'$. Moreover the image $\mathscr{F}\pi_{I'} := \{ (C_{\bm{\delta}})\pi_{I'} \ | \ \bm{\delta} \in X_0 \}$ of $\mathscr{F}$ under this map is a partition of $\PP'$. 

For each $\bm{\delta} \in X_0$, let $\bm{\delta}' := \bm{\delta}\pi_{J'}$, and observe that $(C_{\bm{\delta}})\pi_{I'} = \{ \varepsilon' \in \PP' \ | \ \varepsilon|_{J'} = \bm{\delta}' \} =: C'_{\bm{\delta}'}$. As $\bm{\delta}$ varies over $X_0$ its projection $\bm{\delta}'$ varies over $\E_{J'}$, so $\mathscr{F}\pi_{I'} = \{ C'_{\bm{\delta}'} \ | \ \bm{\delta}' \in \E_{J'} \}$. Let $y_{\bm{\delta}'} := |B' \cap C_{\bm{\delta}'}|$. Then $y_{\bm{\delta}'} = |(B \cap C_{\bm{\delta}})\pi_{I'}|$, and since $\pi_{I'}|_{C_{0_s}}$ is a bijection $C_{0_s} \to \PP'$, we have $|(B \cap C_{\bm{\delta}})\pi_{I'}| = |B \cap C_{\bm{\delta}}|$, so $y_{\bm{\delta}'} = x_{\bm{\delta}}$. Thus
    \begin{equation} \label{eq:x-y}
    \sum_{\bm{\delta} \in X_0} x_{\bm{\delta}}^2 = \sum_{\bm{\delta}' \in \E_{J'}} y_{\bm{\delta}'}^2.
    \end{equation}
Furthermore, by our assumption \eqref{k'}  the structure $\DD{B'}{G'}$ is a $2$-design. Hence by Theorem \ref{t:2des}, the parameters $y_{\bm{\delta}'}$ satisfy the condition \eqref{eq:2des} for each nonempty proper subset $J'$ of $I'$, and also for $J'=\varnothing$ by Lemma~\ref{lem:Jtrivial}. Note that $J' \subsetneq I'$ because $J \subsetneq I$. Setting $k' := |B'|$, $v' := |\PP'|$, and $I' := I \setminus \{s\}$, we then have
    \[
    \sum_{\bm{\delta}' \in \E_{J'}} y_{\bm{\delta}'}^2
    = k' + \frac{k'(k'-1)}{v'-1} \left( \bigg(\prod_{i \in I' \setminus J'} e_i\bigg) - 1 \right).
    \]
Now $k' = p^{2^{s-2}} + 1$ by assumption \eqref{k'}, and $v' = \prod_{i=1}^{s-1} e_i = p^{2^{s-1}} + p^{2^{s-2}} + 1$ by Lemma \ref{lem:e} \eqref{e-prod}. Thus $v'-1=k'(k'-1)$, and substituting into the above equation gives
    \begin{equation} \label{eq:y}
    \sum_{\bm{\delta}' \in \E_{J'}} y_{\bm{\delta}'}^2
    = p^{2^{s-2}} + \prod_{i \in I' \setminus J'} e_i.
    \end{equation}
Note that $I' \setminus J' = (I \setminus \{s\}) \setminus (J \setminus \{s\}) = I \setminus J = J^\c$. Thus Claim 2 follows from \eqref{eq:x-y} and \eqref{eq:y}.

\smallskip
Note that $k(k-1) = \big( p^{2^{s-1}} + 1 \big)p^{2^{s-1}} = p^{2^{s}} + p^{2^{s-1}} = \big( \prod_{i=1}^s e_i \big) - 1 = v-1$,  by Lemma~\ref{lem:e}(a). Then, using \eqref{sum:X+X} and Claims 1 and 2,  we have
    \[
    \sum_{\bm{\delta} \in \E_J} x_{\bm{\delta}}^2
    = \left( p^{2^{s-1}} - p^{2^{s-2}} \right) + \left( p^{2^{s-2}} + \prod_{i \in J^\c} e_i \right)
    = p^{2^{s-1}} + \prod_{i \in J^\c} e_i
    = k + \frac{k(k-1)}{v-1} \left( \prod_{i \in J^\c} e_i - 1 \right).
    \]
Thus equation \eqref{eq:2des} holds for $J$. This completes the proof that \eqref{eq:2des}  holds for any proper subset of $I$ containing $s$. Hence $\DD{B}{G}$ is a $2$-design, completing the proof of Proposition~\ref{lem:B'}.
\end{proof}

\subsection{The case $s=2$}\label{s:s2}

For the construction with $s=2$, we follow a somewhat degenerate form of the approach outlined in Subsection~\ref{s:rec}, in that we take the input design $\DD{B'}{G'}$ to be a complete design on $\PP'=\mathbb{Z}_{e_1}$, where $e_1=p^2+p+1$, with all $(p+1)$-subsets as blocks. We note that Proposition~\ref{lem:B'}\eqref{k'} holds for this design.

\begin{construction} \label{con:s2}
Let $p > 1$ and let $e_1$ and $e_2$ be as in \eqref{def:e}. Let $\PP = \mathbb{Z}_{e_1} \times \mathbb{Z}_{e_2}$ and take $B = B_1 \cup B_2$ to be the set of points with coordinates listed in Table \ref{B:s2-coord}. Thus $B_1=B'\times\{0\}$ with $B'$ a block of the complete design $\DD{B'}{G'}$ mentioned above. Define $\Des{2}{p} = \DD{B}{G}$ as in Definition \ref{def:Des}.
\end{construction}

    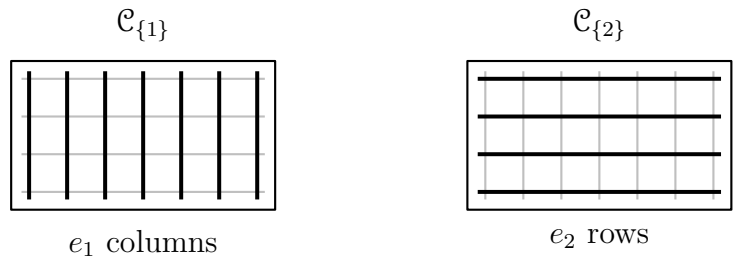
\begin{figure}[ht]
    \begin{center}
    \begin{pspicture}(-4,-2)(4,2)
    \twogridparts
    \end{pspicture}
    \end{center}
    \caption{Partitions $\C_J$ of $\PP = \prod_{i=1}^2 \mathbb{Z}_{e_i}$ for $\varnothing \neq J \subsetneq \{1,2\}$}
    \label{partition:s2}
    \end{figure}

    \begin{table}[ht]
        \centering
        \begin{tabular}{cll}
        \hline
        Subset & Points & Conditions \\
        \hline\hline
        $B_1$ & $(a, 0)$ & $0 \leq a \leq p$ \\
        \hline
        $B_2$ & $(p + a, 2a - 1)$ & $1 \leq a \leq \frac{p^2 - p}{2}$ \\
         & $(p + a, 2a)$ & \\
        \hline
        \end{tabular}
        \caption{Coordinates of points of $B$ in Construction \ref{con:s2}}
        \label{B:s2-coord}
    \end{table}

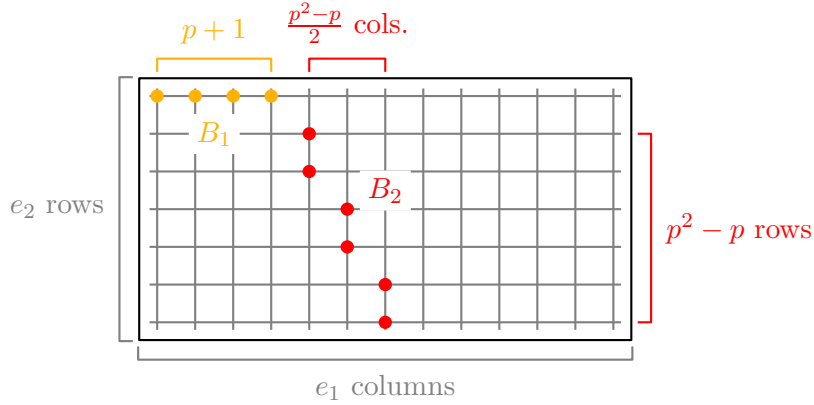
\begin{figure}
        \centering
        \begin{pspicture}(-3,-2.5)(3,2.5)
        \twogridB
        \end{pspicture}
        \caption{Block $B$ for Construction \ref{con:s2}}
        \label{B:s2}
    \end{figure}

We visualise the point set $\PP$ as a grid with $e_1$ columns and $e_2$ rows, where each column is a $\C_{\{1\}}$-part and each row is a $\C_{\{2\}}$-part, and illustrate the block $B$ in Figure \ref{B:s2}. The only proper nonempty subsets of $I$ are the singletons $J = \{i\}$ for $i \in \{1,2\}$, and for each of these $\E_J = \mathbb{Z}_{e_i}$.

The nonzero array values $x_{\bm{\delta}}$, for $\bm{\delta} \in \mathbb{Z}_{e_1} \cup \mathbb{Z}_{e_2}$, are listed in Table \ref{array:s2}. These can be obtained by referring to Table \ref{B:s2-coord} or Figure \ref{B:s2}. Recall that if $\bm{\delta} \in \mathbb{Z}_{e_1}$ (respectively, $\bm{\delta} \in \mathbb{Z}_{e_2}$), the parameter $x_{\bm{\delta}}$ is the number of points of $B$ that lie in the column (respectively, row) $C_{\bm{\delta}}$. We omit from the table those columns and rows which contain no point of $B$.

    \begin{table}[ht]
        \centering
        \begin{tabular}[t]{|crr|}
        \hline
        $\delta_1$ & $\# \delta_1$ & $x_{(\delta_1)}$ \\
        \hline
        $\delta_1 \leq p$ & $p+1$ & $1$ \\
        $p < \delta_1 \leq \frac{p^2 + p}{2}$ & $\frac{p^2 - p}{2}$ & $2$ \\
        \hline
        \end{tabular}
        \quad
        \begin{tabular}[t]{|crr|}
        \hline
        $\delta_2$ & $\# \delta_2$ & $x_{(\delta_2)}$ \\
        \hline
        $0$ & $1$ & $p+1$ \\
        $\delta_2 \geq 1$ & $p^2 - p$ & $1$ \\
        \hline
        \end{tabular}
        \caption{Array values of $B$ in Construction \ref{con:s2}}
        \label{array:s2}
    \end{table}

\begin{remark}\label{r:Des(2,p)}
The design $\Des{2}{p}$ could have been constructed using the combinatorial construction developed in \cite{grids22}: that construction is based on a subgraph $\Delta$ of the complete bipartite  graph $\mathbf{K}_{e_1,e_2}$ isomorphic to $\frac{p^2-p}{2}\mathbf{P}_2 + \mathbf{K}_{1, p+1} + \frac{p^2+p}{2}\mathbf{K}_1$, where the part of the bipartition of size $e_1=p^2+p+1$ consists of the $\frac{p^2 + p}{2}$ isolated vertices, the $\frac{p^2 - p}{2}$ central-vertices of the paths $\mathbf{P}_2$, and the $p+1$ non-central-vertices of the $K_{1,p+1}$; and the part of size $e_2 = p^2 - p + 1$ consists of the $p^2 - p$ end-vertices of the paths $\mathbf{P}_2$ and the central-vertex of the claw $\mathbf{K}_{1, p+1}$. The automorphism group of this graph $\Delta$ is a subgroup of $G=S_{e_1}\times S_{e_2}$, namely 
    \[
    G_\Delta=S_{p+1}\times S_{(p^2+p)/2} \times (S_2\wr S_{(p^2-p)/2}).
    \]
We show in Proposition~\ref{p:Des(2,p)} that $\Des{2}{p}$ is a $2$-design with $v=k(k-1)$; given that result it follows from \cite[Proposition 3.1(a)]{grids22} that the value of $\lambda$ is 
    \[
    \lambda= \frac{(e_1-1)!\cdot (e_2-1)!}{|G_\Delta|} =  \frac{(p^2+p)!\cdot (p^2-p)!}{(p+1)!\cdot ((p^2+p)/2)!\cdot 2^{(p^2-p)/2}((p^2-p)/2)!}. 
    \]
We have not been able to find this family of designs in the literature; it appears to be new.  
\end{remark}

\begin{proposition} \label{p:Des(2,p)}
For each integer $p > 1$, the design $\Des{2}{p}$ in Construction \emph{\ref{con:s2}} is a $2$-grid-imprimitive, block-transitive $2-(p^4+p^2+1, p^2+1, \lambda)$ design relative to $G, \lambda$ as in Remark~\emph{\ref{r:Des(2,p)}}.
\end{proposition}

\begin{proof}
By construction $G$ acts $2$-grid-imprimitively and block-transitive on  $\Des{2}{p}$. By Theorem \ref{t:2des}, to prove that $\Des{2}{p}$ is a $2$-design it suffices to show that conditions \eqref{eq:2des} hold for all proper, nonempty subsets $J$ of $I = \{1,2\}$. The only such subsets are the singletons. We can see from Table \ref{array:s2} that points in the subset $B_1$ all have the same values of $x_{\delta|_J}$ for $J = \{1\}$ or $J=\{2\}$, and likewise for points in $B_2$. The set $B$ is the disjoint union of $B_1$ and $B_2$, so $k:=|B| = |B_1| + |B_2| = (p+1) + (e_2 - 1) = p^2 + 1$. Also $v=|\PP|=e_1e_2= p^4+p+2+1$, so $v=k(k-1)$. Moreover, for any $i \in \{1,2\}$, the quantity $\sum_{(\delta_i) \in \E_{\{i\}}} x_{(\delta_i)}^2 = \sum_{\delta_i \in \mathbb{Z}_{e_i}} x_{(\delta_i)}^2$ can be obtained by taking the sum of the products of the entries in the column `\# $C_{(\delta_i)}$' with the squares of the entries in the column `$x_{(\delta_i)}$' of Table \ref{array:s2}. Thus, using \eqref{def:e} to simplify the resulting sums, we get
    \[
    \sum_{\delta_1 \in \mathbb{Z}_{e_1}} x_{(\delta_1)}^2
    = (p+1) \cdot 1^2 + \left(\frac{e_2 - 1}{2}\right) \cdot 2^2
    = (p+1) + (p^2 - p) \cdot 2
    = k + (e_2 - 1)
    \]
which is condition \eqref{eq:2des} for $J = \{1\}$, and
    \[
    \sum_{\delta_2 \in \mathbb{Z}_{e_2}} x_{(\delta_2)}^2
    = 1 \cdot (p+1)^2 + (e_2 - 1) \cdot 1^2
    = (p^2 + 2p + 1) + (p^2 - p)
    = k + (e_1 - 1),
    \]
which is condition \eqref{eq:2des} for $J = \{2\}$. Therefore, by Theorem \ref{t:2des}, the design $\Des{2}{p}$ is a $2$-design. 
\end{proof}

\subsection{The case $s = 3$}\label{s:s3}

To illustrate the construction with a nondegenerate input design $\DD{B'}{G'}$, we use the designs $\Des{2}{p}$ to construct $3$-grid-imprimitive designs with the approach of Subsections~\ref{s:rec} and~\ref{s:induct}.

\begin{construction} \label{con:s3}
Let $p > 1$ and let $e_1$, $e_2$, and $e_3$ be as in \eqref{def:e}. Let $\PP = \mathbb{Z}_{e_1} \times \mathbb{Z}_{e_2} \times \mathbb{Z}_{e_3}$ and take $B := (B'\times\{0\}) \cup B_3$, where $B'$ is the set of points in Table \ref{B:s2-coord} and $B_3$ is as follows: for $p = 2$,
    \[
    B_3 := \{ (3 + a,b,c-1), (3 + a,b,c)  \ | \ 1 \leq a \leq 3, \ 1 \leq b \leq 2, \ c = 2(a + 3(b - 1)) \}
    \]
and for $p > 2$, $B_3 := \bigcup_{i=1}^4 B_{3,i}$ where the sets $B_{3,i}$ consist of points with coordinates listed in Tables \ref{B:s3odd-coord} (for odd $p$) or \ref{B:s3even-coord} (for even $p > 2$). We abuse notation slightly and, for $B_1$ and $B_2$ the two-dimensional sets given in Table \ref{B:s2-coord}, we also denote the three-dimensional sets $B_1 \times \{0\}$ and $B_2 \times \{0\}$ by $B_1$ and $B_2$, respectively. Define $\Des{3}{p} := \DD{B}{G}$ as in Definition \ref{def:Des}, with $G=S_{e_1}\times S_{e_2}\times S_{e_3}$.
\end{construction}

    \begin{table}[ht]
        \centering
        \begin{tabular}{lll}
        \hline
        Subset & Points & Conditions \\
        \hline\hline
        $B_{3,1}$ & \parbox[t]{3cm}{$(a, b, c-1)$ \\ $(a,b,c)$}
            & $1 \leq a \leq p$, \ $1 \leq b \leq \frac{p^2 - 1}{2}$, \ $c = 2\big( a + (b-1)p \big)$ \\
        \hline
        $B_{3,2}$ & \parbox[t]{3cm}{$(p+a, b+\delta, c-1)$ \\ $(p+a,b+\delta,c)$}
            & \parbox[t]{8.5cm}{$1 \leq a \leq \frac{p^2 - p}{2}$, \ $1 \leq b \leq \frac{p^2 - p}{2}$, \\
            $c = 2p \left( \frac{p^2 - 1}{2} \right) + 2\big( a + (b-1)(p^2-p) \big)$, \\ 
            $\delta=\frac{p^2 - p}{2}$ if $1~\leq~a~\leq~\lfloor \frac{p^2 - p + 2}{4}\rfloor$ and $b \in \{2a-1, \\ 2a\}$, \ and \ $\delta= 0$ \ otherwise} \\
        \hline
        $B_{3,3}$ & \parbox[t]{5cm}{$\left( \frac{p^2 + p}{2} + a, \ \frac{p^2 - 1}{2} + b, \ c-1 \right)$ \\ $\left( \frac{p^2 + p}{2} + a, \ \frac{p^2 - 1}{2} + b, \ c \right)$}
            & \parbox[t]{8.5cm}{$1 \leq a \leq p$, \ $1 \leq b \leq \frac{p^2 - 2p + 1}{2}$, \\
            $c = 2\left(\frac{p^4 + p^2 - 2p}{4}\right) + 2\big( a + (b-1)p \big)$} \\
        \hline
        $B_{3,4}$ & 
        \parbox[t]{5.25cm}{$\left( \frac{p^2 + 3p}{2} + a, \  b+\delta, \ c-1 \right)$ \\ $\left( \frac{p^2 + 3p}{2} + a, \  b+\delta, \ c \right)$}
            & \parbox[t]{8.5cm}{$1 \leq a \leq \frac{p^2 - p}{2}$, \ $1 \leq b \leq \frac{p^2 - p}{2}$, \\
            $c = 2 \left( \frac{p^4 + 2p^3 - 3p^2}{4} \right) + 2\big( a + (b-1)(p^2 - p) \big)$, \\
            $\delta = \frac{p^2 - p}{2} $\ if\ $a=1$ and $\delta=0$ otherwise} \\
        \hline
        \end{tabular}
        \caption{Coordinates of points in $\bigcup_{i=1}^4 B_{3,i}$ in Construction \ref{con:s3} for odd $p$}
        \label{B:s3odd-coord}
    \end{table}

    \begin{table}[ht]
        \centering
        \begin{tabular}{lll}
        \hline
        Subset & Points & Conditions \\
        \hline\hline
        $B_{3,1}$ & \parbox[t]{5cm}{$(a, b, c-1)$ \\ $(a,b,c)$}
            & $1 \leq a \leq p$, \ $1 \leq b \leq \frac{p^2 - 2p}{2}$, \ $c = 2\big( a + (b-1)p \big)$ \\
        \hline
        $B_{3,2}$ & \parbox[t]{5cm}{$(p+a, b+\delta, c-1)$ \\ $(p+a,b+\delta,c)$}
            & \parbox[t]{8.5cm}{$1 \leq a \leq \frac{p^2 - p}{2}$, \ $1 \leq b \leq \frac{p^2 - p}{2}$, \\
            $c = 2p \left( \frac{p^2 - 2p}{2} \right) + 2\big( a + (b-1)(p^2-p) \big)$ \\ 
            $\delta=\frac{p^2 - p}{2}$ if $1 \leq a \leq \lfloor\frac{p^2 - p}{4}\rfloor$ and $b \in \{ 2a-1, 2a \}$, \ and $\delta= 0$ otherwise} \\
        \hline
        $B_{3,3}$ & \parbox[t]{5.25cm}{$\left( \frac{p^2 + p}{2} + a, \ \frac{p^2 - 2p}{2} + b, \ c-1 \right)$ \\ $\left( \frac{p^2 + p}{2} + a, \ \frac{p^2 - 2p}{2} + b, \ c \right)$}
            & \parbox[t]{8.5cm}{$1 \leq a \leq p$, \ $1 \leq b \leq \frac{p^2}{2}$, \\
            $c = 2\left(\frac{p^4 - 3p^2}{4}\right) + 2\big( a + (b-1)p \big)$} \\
        \hline
        $B_{3,4}$ & \parbox[t]{5cm}{$\left( \frac{p^2 + 3p}{2} + a,  b+\delta, c-1\right)$ \\ $\left( \frac{p^2 + 3p}{2} + a, b+\delta,c\right)$}
            & \parbox[t]{8.5cm}{$1 \leq a \leq \frac{p^2 - p}{2}$, \ $1 \leq b \leq \frac{p^2 - p}{2}$, \\ $c = 2 \left( \frac{p^4 + 2p^3 - 3p^2}{4} \right) + 2\big( a + (b-1)(p^2 - p) \big)$, \\ $\delta=\frac{p^2 - p}{2}$\ if \ $a=1$ \ and \ $\delta= 0$ \ otherwise} \\
        \hline
        \end{tabular}
        \caption{Coordinates of points in $\bigcup_{i=1}^4 B_{3,i}$ in Construction \ref{con:s3} for even $p > 2$}
        \label{B:s3even-coord}
    \end{table}

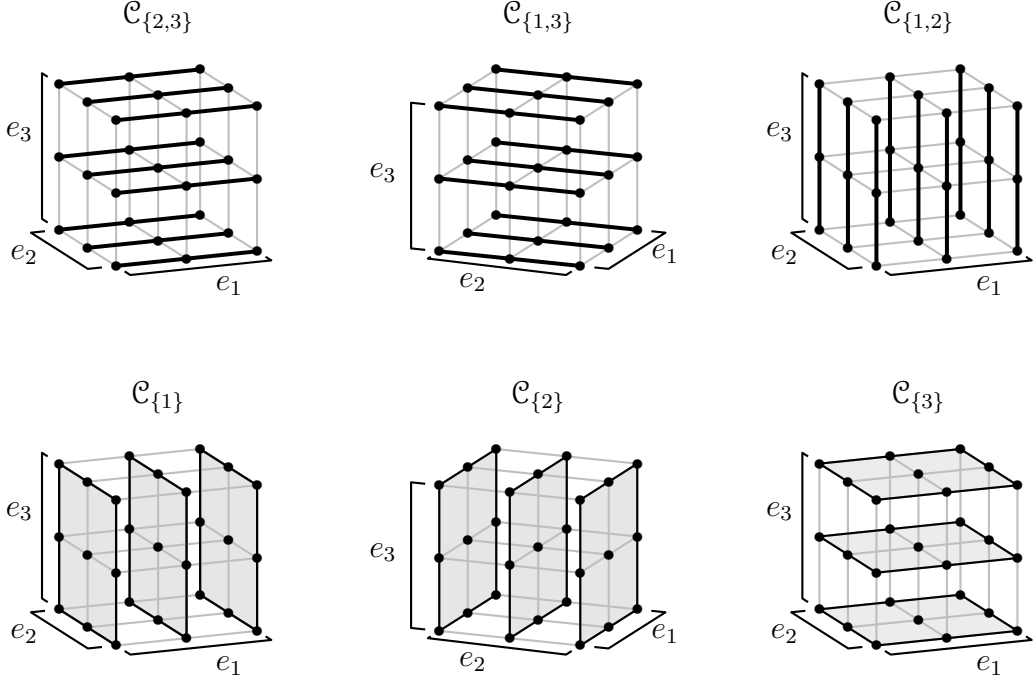
\begin{figure}[ht]
    \centering
    \begin{pspicture}(-4,-2.5)(4,2.5)
    \threegridlines
    \end{pspicture}
    
    \begin{pspicture}(-4,-2.5)(4,2.5)
    \threegridplanes
    \end{pspicture}
    \caption{Partitions $\C_J$ of $\PP = \prod_{i=1}^3 \mathbb{Z}_{e_i}$ for $\varnothing \neq J \subsetneq \{1,2,3\}$}
    \label{partition:s3}
    \end{figure}

We visualise $\PP$ as a three-dimensional grid, with $e_3$ ``layers'' labelled by elements of $\mathbb{Z}_{e_3}$, each of which is a two-dimensional grid with $e_1$ columns labelled by $\mathbb{Z}_{e_1}$ and $e_2$ rows labelled by $\mathbb{Z}_{e_2}$. A $\C_{\{3\}}$-part $C_{(\delta_3)}$ is the set of all points in the same layer $\delta_3$, a $\C_{\{1\}}$-part $C_{(\delta_1)}$ is the set of all points in the same column $\delta_1$ and any layer, and a $\C_{\{2\}}$-part is the set of all points in the same row $\delta_2$ and any layer. A $\C_{\{1,3\}}$-part $C_{(\delta_1,\delta_3)}$ is the set of all points in column $\delta_1$ and layer $\delta_3$, a $\C_{\{2,3\}}$-part $C_{(\delta_2,\delta_3)}$ is the set of all points in row $\delta_2$ and layer $\delta_3$, and a $\C_{\{1,2\}}$-part $C_{(\delta_1,\delta_2)}$ is the set of all points in column $\delta_1$ and row $\delta_2$ (and any layer). These partitions are illustrated in Figure \ref{partition:s3}.

The generating block $B$ is illustrated in Figure \ref{B:s3p2} (for $p=2$), Figure \ref{B:s3odd} (for odd $p$), and Figure \ref{B:s3even} (for even $p > 2$). To save space, each of these figures shows the projection of $B$ onto a single two-dimensional grid representing $\mathbb{Z}_{e_1} \times \mathbb{Z}_{e_2}$. In Figure \ref{B:s3p2} the five points in $B'\times\{0\} = B_1 \cup B_2$ all lie in the layer $\delta_3 = 0$. Each dot in $B_3$ represents two points with the same $(\delta_1,\delta_2)$-coordinates but which lie in different $\delta_3$-layers, and each layer $\delta_3 \in \{1, \ldots, 12\}$ contains exactly one point of $B_3$. In Figures \ref{B:s3odd} and \ref{B:s3even}, the subsets of $B$ are represented as regions, which indicate a point in every intersection of a row and a column inside that region. (For instance, the region $B_{1,2}$ consists of one row that intersects $p$ columns, so the subset $B_{1,2}$ contains $p$ points that lie in the same row in consecutive columns.) The set $B_1$ is divided into two smaller subsets $B_{1,1}$ and $B_{1,2}$. Again the points in $B' \times \{0\}$ lie in the layer $\delta_3 = 0$. Each point $(\delta_1,\delta_2)$ in the regions corresponding to $B_3 = \bigcup_{i=1}^4 B_{3,i}$ represents two elements $(\delta_1,\delta_2,\delta_3)$ of $B_3$ that lie in distinct $\delta_3$-layers, and each layer $\delta_3 \in \{1, \ldots, e_3 - 1\}$ contains exactly one point of $B_3$.

    \begin{figure}
    \centering
    \begin{pspicture}(-2,-1.5)(2,1.5) \psset{unit=1.25}
    \threeblock \threegridB
    \end{pspicture}
    \caption{Block $B$ for Construction \ref{con:s3} with $p=2$, projected onto $\mathbb{Z}_7 \times \mathbb{Z}_3$}
    \label{B:s3p2}
    \end{figure}
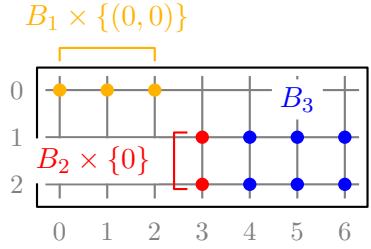

    \begin{figure}[ht]
        \centering
        \begin{subfigure}{0.8\textwidth}
            \begin{pspicture}(-6,-5)(6,5)
            \threegridBgenImodIV
            \end{pspicture}
        \caption{Case $p \equiv 1 \pmod{4}$}
        \label{B:s3p1mod4}
        \end{subfigure}
        
        \begin{subfigure}{0.8\textwidth}
        \centering
            \begin{pspicture}(-6,-6)(6,5)
            \threegridBgenImodIV
            \end{pspicture}
        \caption{Case $p \equiv 3 \pmod{4}$}
        \label{B:s3p3mod4}
        \end{subfigure}
        
        \caption{Block $B$ for Construction \ref{con:s3} with odd $p$, projected onto $\mathbb{Z}_{e_1} \times \mathbb{Z}_{e_2}$}
        \label{B:s3odd}
    \end{figure}

    \begin{figure}[ht]
        \centering
        \begin{subfigure}{0.8\textwidth}
        \centering
            \begin{pspicture}(-6,-5)(6,5)
            \threegridBgenOmodIV
            \end{pspicture}
        \caption{Case $p \equiv 0 \pmod{4}$}
        \label{B:p=0mod4}
        \end{subfigure}
        
        \begin{subfigure}{0.8\textwidth}
        \centering
            \begin{pspicture}(-6,-6)(6,5)
            \threegridBgenIImodIV
            \end{pspicture}
        \caption{Case $p \equiv 2 \pmod{4}$}
        \label{B:p=2mod4}
        \end{subfigure}
        
        \caption{Block $B$ for Construction \ref{con:s3} with even $p > 2$, projected onto $\mathbb{Z}_{e_1} \times \mathbb{Z}_{e_2}$}
        \label{B:s3even}
    \end{figure}

In view of Proposition \ref{lem:B'}, to show that $\Des{3}{p}$ is a $2$-design we only need to verify conditions \eqref{eq:2des} for the subsets $J = \{1\}$, $\{2\}$, and $\{1,2\}$. First we determine the array values of $B$ for these subsets.

\begin{lemma}\label{lem:s3array}
    The array values of the subset $B$ in Construction~$\ref{con:s3}$ corresponding to  the subsets $J = \{1\}$, $\{2\}$, and $\{1,2\}$, are  listed in Table \ref{array:s3p2} for $p = 2$, Table \ref{array:s3odd} for odd $p$, and Table \ref{array:s3even} for even $p > 2$ (where in Tables~\ref{array:s3odd} and~\ref{array:s3even},  $B_{1,1}=\{(0,0,0)\}$ and $B_{1,2}=\{(a,0,0)\mid 1\leq a\leq p\}$).
\end{lemma}

    \begin{table}[ht]
        \centering
        \begin{tabular}[t]{|crr|}
        \hline
        Subset & $\# \delta|_{\{1,2\}}$ & $x_{\delta|_{\{1,2\}}}$ \\
        \hline
        $B' \times \{0\}$ & $5$ & $1$ \\
        $B_3$ & $6$ & $2$ \\
        \hline
        \end{tabular}
        \quad
        \begin{tabular}[t]{|crr|}
        \hline
        $\delta_2$ & $\# \delta_2$ & $x_{\delta_2}$ \\
        \hline
        $0$ & $1$ & $3$ \\
        $1,2$ & $2$ & $7$ \\
        \hline
        \end{tabular}
        \quad
        \begin{tabular}[t]{|crr|}
        \hline
        $\delta_1$ & $\# \delta_1$ & $x_{\delta_1}$ \\
        \hline
        $0,1,2$ & $3$ & $1$ \\
        $3$ & $1$ & $2$ \\
        $4,5,6$ & $3$ & $4$ \\
        \hline
        \end{tabular}
        \caption{Array values for $B$ in Construction \ref{con:s3} with $p=2$}
        \label{array:s3p2}
    \end{table}
    
    \begin{table}[ht]
        \centering
        \begin{tabular}[t]{|crr|}
        \hline
        $\bm{\delta} = \delta|_{\{1,2\}}$ & $\# \bm{\delta}$ & $x_{\bm{\delta}}$ \\
        \hline
        $\delta \in B' \times \{0\}$ & $p^2 + 1$ & $1$ \\
        $\delta \in B_3$ & $\frac{e_3 - 1}{2}$ & $2$ \\
        \hline
        \end{tabular}
        \begin{tabular}[t]{|crr|}
        \hline
        $\delta_2$ & $\# \delta_2$ & $x_{\delta_2}$ \\
        \hline
        $0$ & $1$ & $p+1$ \\
        $\delta_2 \geq 1$ & $p^2 - p$ & $e_1$ \\
        \hline
        \end{tabular}
        \begin{tabular}[t]{|crr|}
        \hline
        $\delta_1$ & $\# \delta_1$ & $x_{(\delta_1)}$ \\
        \hline
        $0$ & $1$ & $1$ \\
        $0 < \delta_1 \leq p$ & $p$ & $p^2$ \\
        $p < \delta_1 \leq \frac{p^2 + p}{2}$ & $\frac{p^2 - p}{2}$ & $e_2 + 1$ \\
        $\frac{p^2 + p}{2} < \delta_1 \leq \frac{p^2 + 3p}{2}$ & $p$ & $(p-1)^2$ \\
        $\frac{p^2 + 3p}{2} < \delta_1 \leq p^2 + p$ & $\frac{p^2 - p}{2}$ & $e_2 - 1$ \\
        \hline
        \end{tabular}
        \caption{Array values for $B$ in Construction \ref{con:s3} with odd $p$}
        \label{array:s3odd}
    \end{table}

    \begin{table}[ht]
        \centering
        \begin{tabular}[t]{|crr|}
        \hline
        $\bm{\delta} = \delta|_{\{1,2\}}$ & $\# \bm{\delta}$ & $x_{\bm{\delta}}$ \\
        \hline
        $\delta \in B' \times \{0\}$ & $p^2 + 1$ & $1$ \\
        $\delta \in B_3$ & $\frac{e_3 - 1}{2}$ & $2$ \\
        \hline
        \end{tabular}
        \begin{tabular}[t]{|crr|}
        \hline
        $\delta_2$ & $\# \delta_2$ & $x_{\delta_2}$ \\
        \hline
        $0$ & $1$ & $p+1$ \\
        $\delta_2 \geq 1$ & $p^2 - p$ & $e_1$ \\
        \hline
        \end{tabular}
        \begin{tabular}[t]{|crr|}
        \hline
        $\delta_1$ & $\# \delta_1$ & $x_{\delta_1}$ \\
        \hline
        $0$ & $1$ & $1$ \\
        $1 \leq \delta_1 \leq p$ & $p$ & $(p-1)^2$ \\
        $p < \delta_1 \leq \frac{p^2 + p}{2}$ & $\frac{p^2 - p}{2}$ & $e_2 + 1$ \\
        $\frac{p^2 + p}{2} < \delta_1 \leq \frac{p^2 + 3p}{2}$ & $p$ & $p^2$ \\
        $\frac{p^2 + 3p}{2} < \delta_1 \leq p^2 + p$ & $\frac{p^2 - p}{2}$ & $e_2 - 1$ \\
        \hline
        \end{tabular}
        \caption{Array values for $B$ in Construction \ref{con:s3} with even $p > 2$}
        \label{array:s3even}
    \end{table}

\begin{proof}
These values are obtained by referring to Construction \ref{con:s3}, or they may be seen visually by referring to Figures \ref{B:s3p2}, \ref{B:s3odd}, and \ref{B:s3even}, as well as Tables~\ref{B:s2-coord}, \ref{B:s3odd-coord} and~\ref{B:s3even-coord}.

(a) The parameter $x_{(\delta_1,\delta_2)}$, for each $(\delta_1,\delta_2) \in \mathbb{Z}_{e_1} \times \mathbb{Z}_{e_2}$, is the number of points of $B$ that belong to row $\delta_1$ and column $\delta_2$ and to any layer, so in the figures this is the number of points of $B$ represented by each dot $(\delta_1,\delta_2)$. For $(\delta_1,\delta_2) \in B'$, this number is $1$, and for $(\delta_1,\delta_2) \in (B_3)\pi_{\{1,2\}} = \bigcup_{i=1}^4 (B_{3,i})\pi_{\{1,2\}}$, this number is $2$.

(b) The parameter $x_{(\delta_2)}$ is the number of points of $B$ that belong to the column $\delta_2$ and to any layer, so in the figures $x_{(\delta_2)}$ is the number of points of $B$ that belong to row $\delta_2$. The top row $\delta_2 = 0$ contains all points $B_1 \times \{(0,0)\}$ and nothing else, so for $\delta_2 = 0$, we have $x_{(\delta_2)} = |B_1| = p+1$ (likewise $x_{(\delta_2)} = 3 = 2 + 1$ for $p = 2$). Each row $\delta_2 \in \{1, \ldots, e_2-1\}$ contains exactly one point of $B_2\times\{0\}$ and a certain number of points of $B_3$. If $p = 2$, then there are exactly $2 \cdot 3 = 6 = e_1 - 1$ points of $B_3$ in each of these rows, so in this case $x_{(\delta_2)} = 7 = e_1$ for $\delta_2 \in \{1,2\}$. If $p > 2$ then, by adding the widths of the regions $B_{3,i}$ in Figures \ref{B:s3odd} and \ref{B:s3even} that intersect any row $\delta_2 > 0$, we find that the number of points of $B_3$ is
    \[
    2 \cdot \left( p + \left( \frac{p^2-p}{2} - 1 \right) + 1 \right)
    = 2 \cdot \frac{p^2 + p}{2}
    = p^2 + p
    = e_1 - 1.
    \]
Hence, if $p > 2$, then $x_{(\delta_2)} = 1 + (e_1 - 1) = e_1$ for any $\delta_2 \in \{1, \ldots, e_2 - 1\}$.

(c) The parameter $x_{(\delta_1)}$ is the number of points of $B$ that belong to the column $\delta_1$ and to any layer, so in the figures $x_{(\delta_1)}$ is the number of points of $B$ that belong to the column $\delta_1$.

\smallskip\noindent
\emph{The  case $p = 2$:}\quad We see in Figure \ref{B:s3p2} that each of the columns $0,1,2$ contains exactly one point of $B_1$, column $3$ contains the two points of $B_2$, and each of the columns $4,5,6$ contains four points of $B_3$. This gives the values in Table \ref{array:s3p2}.

\smallskip\noindent
\emph{The  case $p > 2$:} \quad
We see in Figures \ref{B:s3odd} and \ref{B:s3even} that the column $0$ contains only the unique point of $B_{1,1}$, so $x_{(\delta_1)} = 1$ for $\delta_1 = 0$. Each of the next $p$ columns contains one point of $B_{1,2}$ and $2m$ points of $B_{3,1}$, where (see Tables \ref{B:s3odd-coord} and \ref{B:s3even-coord})
    \[
    m = \begin{cases}
        \frac{p^2 - 1}{2} &\text{if $p$ is odd} \\
        \frac{p^2 - 2p}{2} &\text{if $p$ is even}.
    \end{cases}
    \]
So for $\delta_1 \in \{1, \ldots, p\}$ we have $x_{(\delta_1)} = 2m+1 = p^2$ or $(p-1)^2$, for $p$ odd or even, respectively. Each of the next $(p^2-p)/2$ columns contains two points of $B_2$ and $2 \cdot (p^2-p)/2 = p^2 - p$ points of $B_{3,2}$, so $x_{(\delta_2)} = p^2 - p + 2 = e_2 + 1$ for $\delta_2 \in \left\{ p+1, \ldots, (p^2 + p)/2 \right\}$. Each of the next $p$ columns contains $2n$ points of $B_{3,2}$, where (see Tables \ref{B:s3odd-coord} and \ref{B:s3even-coord})
    \[
    n = \begin{cases}
        \frac{p^2 - 2p + 1}{2} = \frac{(p-1)^2}{2} &\text{if $p$ is odd} \\
        \frac{p^2}{2} &\text{if $p$ is even}.
    \end{cases}
    \]
So for $\delta_1 \in \left\{ (p^2 + p)/2 + 1, \ldots, (p^2 + 3p)/2 \right\}$ we have $x_{(\delta_1)} = 2n = (p-1)^2$ or $p^2$, if $p$ is odd or even, respectively. Each of the remaining $(p^2 - p)/2$ columns contains $2 \cdot (p^2 - p)/2 = p^2 - p$ points of $B_{3,4}$, so $x_{(\delta_1)} = p^2 - p = e_2 - 1$ 
for $\delta_1 \in \left\{ (p^2 + 3p)/2 + 1, \ldots, p^2 + p \right\}$. This completes the proof.
\end{proof}

\begin{proposition} \label{p:Des(3,p)}
For each integer $p > 1$ and $G$ as in Construction \emph{\ref{con:s3}}, the design $\Des{3}{p}$ as in Construction \emph{\ref{con:s3}} is a $3$-grid-imprimitive, block-transitive  $2-(v,k,\lambda)$ design with $v=p^8+p^4+1$ and $k=p^4+1$, for some $\lambda$.
\end{proposition}

\begin{proof}
The values for $v=|\PP|$ and $k=|B|$ follow from Construction~\ref{con:s3} and Lemma~\ref{lem:e}, and we note that $k(k-1)=v-1$. Also $G$ acts as a $3$-grid-imprimitive, block-transitive group of automorphisms.
We see from the definition of the block $B=(B' \times \{0\})\cup B_3$ in Construction~\ref{con:s3} that $B'$ is the generating block in Construction \ref{con:s2} with the coordinates listed in Table \ref{B:s2-coord}, and $\DD{B'}{G'}=\Des{2}{p}$ is a $2$-design by Proposition~\ref{p:Des(2,p)}. Thus $B$ satisfies condition \eqref{k'} of Proposition \ref{lem:B'}. Moreover since each $\delta_3$-layer, with $\delta_3 \in\mathbb{Z}_{e_3}\setminus\{0\}$, contains a unique point of $B_3$, the set $B_3$ satisfies condition \eqref{B_s} of Proposition \ref{lem:B'}. Therefore by Proposition \ref{lem:B'}, in order to show that $\Des{3}{p}$ is a $2$-design, we only need to verify that condition \eqref{eq:2des} holds for all nonempty subsets $J$ of $\{1,2,3\}$ that do not contain $3$, namely, for $J = \{1\}$, $\{2\}$, and $\{1,2\}$.

\smallskip\noindent\emph{Case $J = \{2\}$:}\quad 
In this case, the right hand side of \eqref{eq:2des} is $k+e_1e_3-1$ since $k(k-1)=v-1$.
By Lemma \ref{lem:s3array} (Tables \ref{array:s3p2}, \ref{array:s3odd}, and \ref{array:s3even}), for any value of $p$ (including $p = 2$), the $p^2 - p$ rows $\delta_2 > 0$ each contain the same number of points of $B$. Thus the sum $\sum_{\delta_2 \in \mathbb{Z}_{e_2}} x_{(\delta_2)}^2$ can be written as
    \[
    \sum_{\delta_2 \in \mathbb{Z}_{e_2}} x_{(\delta_2)}^2
    = \sum_{\delta_2 = 0} x_{(\delta_2)}^2 + \sum_{\delta_2 > 0} x_{(\delta_2)}^2
    = 1 \cdot (p+1)^2 + (p^2 - p) \cdot e_1^2.
    \]
Since $e_1-1=p^2+p$ and $e_2-1=p^2-p$, we have $(e_1 - 1)(e_2 - 1) = p^4-p^2= e_3 - 1$, so
    \[
    (p^2 - p) \cdot e_1^2
    = (e_2 - 1)(e_1 - 1 + 1)e_1
    = (e_3 - 1)e_1 + (e_2 - 1)e_1.
    \]
Also $(p+1)^2 = p^2 + 2p + 1 = e_1 + p$ and $(e_2 - 1)e_1 = (p^2 - p)(p^2 + p + 1) = p^4 - p$, so
    \begin{align*}
    \sum_{\delta_2 \in \mathbb{Z}_{e_2}} x_{(\delta_2)}^2
    &= (p + 1)^2 + (e_3 - 1)e_1 + (e_2 - 1)e_1 = (e_1+p)+(e_1e_3 -e_1)+(p^4-p)\\
    &= p^4 + e_1e_3
    = k + e_1e_3 - 1.
    \end{align*}
Thus condition \eqref{eq:2des} holds for $J = \{2\}$.

\smallskip\noindent\emph{Case $J = \{1,2\}$:}\quad 
Since, for $\bm{\delta} = (\delta_1,\delta_2) \in \mathbb{Z}_{e_1} \times \mathbb{Z}_{e_2}$, the parameter $x_{\bm{\delta}} = 0$ whenever $\bm{\delta} \notin B\pi_{\{1,2\}}$, the sum $\sum_{\bm{\delta} \in \mathbb{Z}_{e_1} \times \mathbb{Z}_{e_2}} x_{\bm{\delta}}^2 = \sum_{\bm{\delta} \in B\pi_{\{1,2\}}} x_{\bm{\delta}}^2$. Since the image $B\pi_{\{1,2\}}$ is the disjoint union of $B'$ and $(B_3)\pi_{\{1,2\}}$, and since by Lemma \ref{lem:s3array} (Tables \ref{array:s3p2}, \ref{array:s3odd}, and \ref{array:s3even}) distinct points in either of these subsets have the same array values $x_{\bm{\delta}}$, the sum $\sum_{\bm{\delta} \in \mathbb{Z}_{e_1} \times \mathbb{Z}_{e_2}} x_{\bm{\delta}}^2$ can be written as
    \begin{align*}
    \sum_{\bm{\delta} \in \mathbb{Z}_{e_1} \times \mathbb{Z}_{e_2}} x_{\bm{\delta}}^2
    &= |B'| \cdot \left( x_{\bm{\delta}}^2 \ \text{for $\bm{\delta} \in B'$} \right)
    + |(B_3)\pi_{\{1,2\}}| \cdot \left( x_{\bm{\delta}}^2 \ \text{for $\bm{\delta} \in (B_3)\pi_{\{1,2\}}$} \right).
    \end{align*}
Note that $k = p^4 + 1$, that $|B'| = |B_1| + |B_2| = (p+1) + (p^2 - p) = p^2 + 1$, and that $|B_3\pi_{\{1,2\}}| = \frac{|B_3|}{2} = \frac{e_3 - 1}{2} = \frac{p^4 - p^2}{2}$ (see Construction \ref{con:s3} and Lemma \ref{lem:e}). Thus, using the array values listed in the tables, we obtain (for all $p$ including $p=2$)
    \[
    \sum_{\bm{\delta} \in \mathbb{Z}_{e_1} \times \mathbb{Z}_{e_2}} x_{\bm{\delta}}^2
    = (p^2 + 1) \cdot 1^2 + \left(\frac{e_3 - 1}{2}\right) \cdot 2^2
    = p^4 + 1 + e_3 - 1
    = k + e_3 - 1.
    \]
Therefore condition \eqref{eq:2des} holds for $J = \{1,2\}$.

\smallskip\noindent\emph{Case $J = \{1\}$:}\quad 
 When $p = 2$, we have $k=2^4+1=17$ and $e_2e_3=3\cdot 13 = 39$, and by Lemma \ref{lem:s3array} (Table \ref{array:s3p2}) the parameter $x_{(\delta_1)}$ is constant over the columns $\delta_1\in\{ 0,1,2\}$ and over the columns $\delta_1 \in\{4,5,6\}$. Thus
    \begin{align*}
    \sum_{\delta_1 \in \mathbb{Z}_7} x_{(\delta_1)}^2
    &= 3 \cdot 1^2 + 1 \cdot 2^2 + 3 \cdot 4^2
    = 3 + 4 + 48 = 55 = k + e_2e_3 - 1, 
   \end{align*}
so condition \eqref{eq:2des} for $J = \{1\}$ is satisfied. If $p > 2$, then by Lemma \ref{lem:s3array} (Tables \ref{array:s3odd} and \ref{array:s3even}) the columns $\delta_1 \in\mathbb{Z}_{e_1}$ can be partitioned into five subsets, listed in the first column of these tables, such that distinct columns in each of these subsets have the same array value $x_{(\delta_1)}$. Observe that for any $p > 2$, $x_{(\delta_1)}^2 =1$ if $\delta_1=0$,
    \[
    \sum_{0 < \delta_1 \leq p} x_{(\delta_1)}^2 + \sum_{\frac{p^2 + p}{2} < \delta_1 \leq \frac{p^2 + 3p}{2}} x_{(\delta_1)}^2
    = p \cdot \left(p^2\right)^2 + p \cdot \left((p-1)^2\right)^2
    = p^5 + p(p-1)^4
    \]
and
    \begin{align*}
    \sum_{p < \delta_1 \leq \frac{p^2 + p}{2}} x_{(\delta_1)}^2 + \sum_{\frac{p^2 + 3p}{2} < \delta_1 \leq p^2 + p} x_{(\delta_1)}^2
    &= \frac{p^2 - p}{2} \cdot (e_2 + 1)^2 + \frac{p^2 - p}{2} \cdot (e_2 - 1)^2 \\
    &= (e_2 - 1)(e_2^2 + 1).
    \end{align*}
Therefore
    \begin{align*}
    \sum_{\delta_1 \in \mathbb{Z}_{e_1}} x_{(\delta_1)}^2
    &= 1 + p^5 + p(p-1)^4 + (e_2 - 1)(e_2^2 + 1).
    \end{align*}
So, to complete the proof, it remains to show, for any $p>2$, that
    \begin{equation}\label{e:des3p}
        1 + p^5 + p(p-1)^4 + (e_2 - 1)(e_2^2 + 1) = k + e_2e_3 - 1 = p^4 + e_2e_3.
    \end{equation}  
The right hand side of \eqref{e:des3p} satisfies
    \[
    {\rm RHS} = p^4 + (p^2-p+1)(p^4-p^2+1) = p^6-p^5+ p^4 +p^3 -p+1. 
    \]
To evaluate the left hand side  of \eqref{e:des3p}, we first note that
    \[
    e_2^2 + 1 = (e_2 - 1)^2 + 2e_2 = (p^2-p)^2 +2(p^2-p+1) = p^4 -2p^3 +3p^2 -2p+2,
    \]
and hence
    \begin{align*}
    {\rm LHS} &= 1+p^5 + (p^2-p)((p-1)^3 + p^4 -2p^3 +3p^2 -2p +2)\\
        &= 1+p^5 + (p^2-p)( p^4 -p^3 +p  +1) = 1 +p^6 -p^5 +p^4 +p^3 -p. 
    \end{align*}
Thus \eqref{e:des3p} holds also when $p>2$. 

\smallskip
We conclude that, for any integer $p > 1$, 
the incidence structure $\Des{3}{p}$ is a $2$-design.
\end{proof}

\subsection{The case $s=4$}\label{s:s4}

When $s=4$, we were only able to carry out this construction method for $p=2$: we use the design $\Des{3}{2}$ to construct a $4$-grid-imprimitive $2$-design.

\begin{example} \label{ex:s4p2}    
Let $p = 2$ and let $e_1$, $e_2$, $e_3$, and $e_4$ be as in \eqref{def:e}. Then $e_1 = 7$, $e_2 = 3$, $e_3 = 13$, and $e_4 = 241$. Let $\PP = \prod_{i=1}^4 \mathbb{Z}_{e_i}$ and take $B = (B' \times \{0\}) \cup B_4$, where $B'$ is the set of points in Figure \ref{B:s3p2} and $B_4 = \bigcup_{i=1}^6 B_{4,i}$ is the set of points listed in Table \ref{B:s4p2-coord}. Again we abuse notation and, for $i \in \{1,2,3\}$, we denote  by the symbol $B_i$ the four-dimensional subset $B_i \times \{0\}$, where $B_i$ is the three-dimensional subset in Figure \ref{B:s3p2}. Define $\Des{4}{2} := \D(B,K)$ as in Definition \ref{def:Des}.
\end{example}

\begin{table}[ht]
    \centering
    \begin{tabular}{llll}
    \hline
    Subset & Coordinates & \multicolumn{2}{l}{Conditions on $c$ and $d$} \\
    \hline\hline
    $B_{4,1}$ 
    & $(0,0,c,d-1)$, $(0,0,c,d)$    & $c \in \{9, \ldots, 12\}$,    & $d = 2(c-8)$  \\ 
    & $(1,0,c,d-1)$, $(1,0,c,d)$    & $c \in \{1, \ldots, 4\}$,     & $d = 2c + 8$  \\ 
    & $(2,0,c,d-1)$, $(2,0,c,d)$    & $c \in \{5, \ldots, 8\}$,     & $d = 2c + 8$  \\ 
    \hline
    $B_{4,2}$ 
    & $(0,1,c,d-1)$, $(0,1,c,d)$    & $c \in \{5, \ldots, 8\}$,     & $d = 2c + 16$ \\ 
    &                               & $c \in \{11,12\}$,            & $d = 2c + 12$ \\ 
    & $(1,1,c,d-1)$, $(1,1,c,d)$    & $c \in \{3,4\}$,              & $d = 2c + 32$ \\ 
    &                               & $c \in \{9, \ldots, 12\}$,    & $d = 2c + 24$ \\ 
    & $(2,1,c,d-1)$, $(2,1,c,d)$    & $c \in \{1, \ldots, 4\}$,     & $d = 2c + 48$ \\ 
    &                               & $c \in \{7,8\}$,              & $d = 2c + 44$ \\ 
    & $(0,2,c,d-1)$, $(0,2,c,d)$    & $c \in \{5, \ldots, 10\}$,    & $d = 2c + 52$ \\ 
    & $(1,2,c,d-1)$, $(1,2,c,d)$    & $c \in \{1,2\}$,              & $d = 2c + 72$ \\ 
    &                               & $c \in \{9, \ldots, 12\}$,    & $d = 2c + 60$ \\ 
    & $(2,2,c,d-1)$, $(2,2,c,d)$    & $c \in \{1, \ldots, 6\}$,     & $d = 2c + 84$ \\ 
    \hline
    $B_{4,3}$ 
    & $(3,0,c,d-1)$, $(3,0,c,d)$    & $c \in \{1, \ldots, 12\}$,    & $d = 2c + 96$  \\ 
    \hline
    $B_{4,4}$ 
    & $(3,1,c,d-1)$, $(3,1,c,d)$    & $c \in \{3,4\}$,              & $d = 2c + 116$ \\ 
    &                               & $c \in \{7,8\}$,              & $d = 2c + 112$ \\ 
    &                               & $c \in \{11,12\}$,            & $d = 2c + 108$ \\ 
    & $(3,2,c,d-1)$, $(3,2,c,d)$    & $c \in \{1,2\}$,              & $d = 2c + 132$ \\ 
    &                               & $c \in \{5,6\}$,              & $d = 2c + 128$ \\ 
    &                               & $c \in \{9,10\}$,             & $d = 2c + 124$ \\ 
    \hline
    $B_{4,5}$
    & $(4,0,c,d-1)$, $(4,0,c,d)$    & $c \in \{9, \ldots, 12\}$,    & $d = 2c + 128$ \\ 
    & $(5,0,c,d-1)$, $(5,0,c,d)$    & $c \in \{1, \ldots, 4\}$,     & $d = 2c + 152$ \\ 
    & $(6,0,c,d-1)$, $(6,0,c,d)$    & $c \in \{5, \ldots, 8\}$,     & $d = 2c + 152$ \\ 
    \hline
    $B_{4,6}$
    & $(4,1,c,d-1)$, $(4,1,c,d)$    & $c \in \{5, \ldots, 8\}$,     & $d = 2c + 160$ \\ 
    &                               & $c \in \{11,12\}$,            & $d = 2c + 156$ \\ 
    & $(5,1,c,d-1)$, $(5,1,c,d)$    & $c \in \{1,2\}$,              & $d = 2c + 180$ \\ 
    &                               & $c \in \{9, \ldots, 12\}$,    & $d = 2c + 168$ \\ 
    & $(6,1,c,d-1)$, $(6,1,c,d)$    & $c \in \{1, \ldots, 4\}$,     & $d = 2c + 192$ \\ 
    &                               & $c \in \{7,8\}$,              & $d = 2c + 188$ \\ 
    & $(4,2,c,d-1)$, $(4,2,c,d)$    & $c \in \{5, \ldots, 10\}$,    & $d = 2c + 196$ \\ 
    & $(5,2,c,d-1)$, $(5,2,c,d)$    & $c \in \{1,2\}$,              & $d = 2c + 216$ \\ 
    &                               & $c \in \{9, \ldots, 12\}$,    & $d = 2c + 212$ \\ 
    & $(6,2,c,d-1)$, $(6,2,c,d)$    & $c \in \{1, \ldots, 6\}$,     & $d = 2c + 228$ \\ 
    \hline
    \end{tabular}
    \caption{Coordinates of points in $\bigcup_{i=1}^6 B_{4,i}$ in Example \ref{ex:s4p2}}
    \label{B:s4p2-coord}
\end{table}

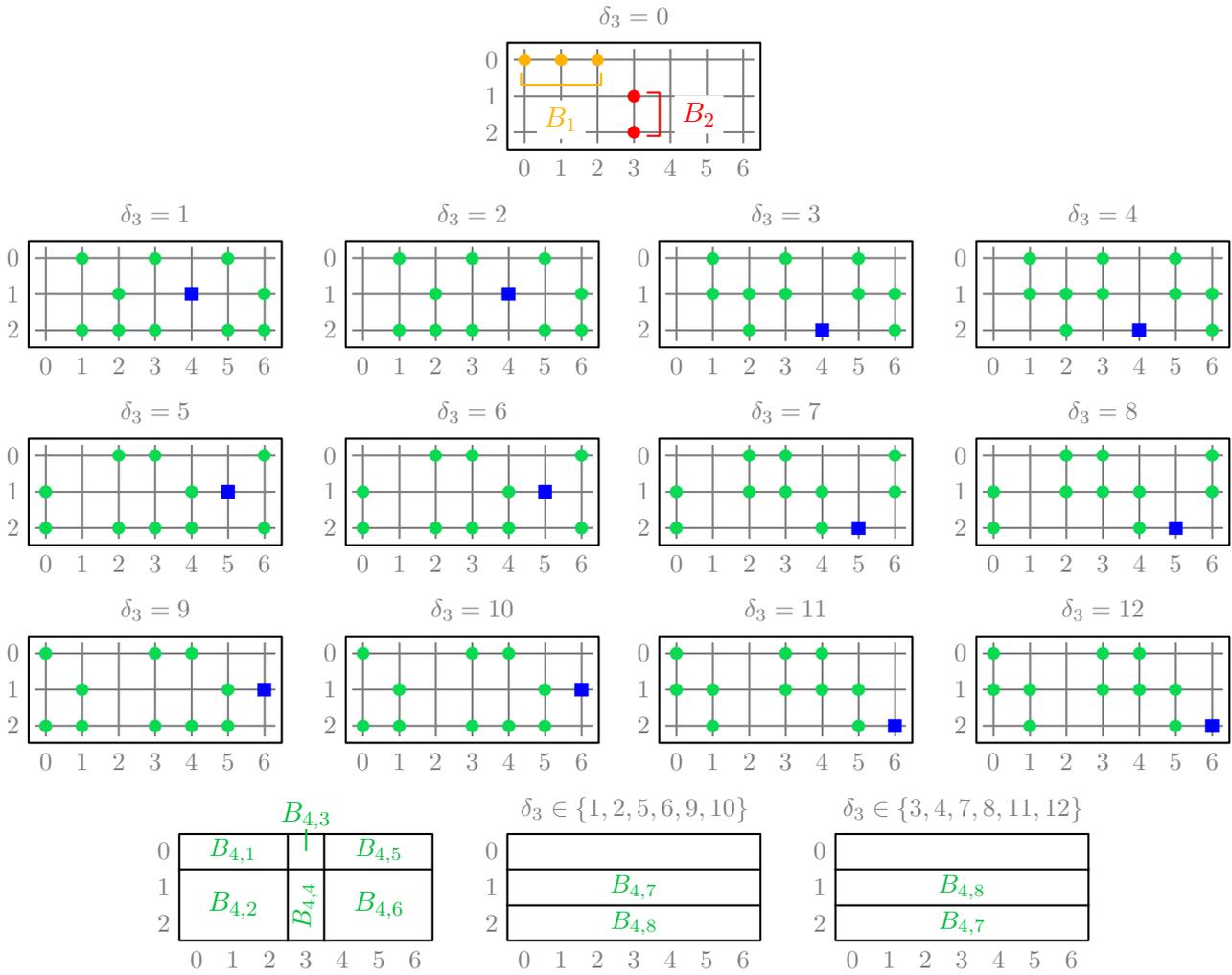
\begin{figure}[ht]
    \begin{center}
    \begin{pspicture}(-6,-1.35)(6,1.35)
    \fourgridBI
    \end{pspicture}
    
    \begin{pspicture}(-6,-1.35)(6,1.35)
    \fourgridBII
    \end{pspicture}
    
    \begin{pspicture}(-6,-1.35)(6,1.35)
    \fourgridBIII
    \end{pspicture}
    
    \begin{pspicture}(-6,-1.35)(6,1.35)
    \fourgridBIV
    \end{pspicture}
    \end{center}
    \caption{Block $B$ for Example \ref{ex:s4p2}, projected onto $\mathbb{Z}_{7} \times \mathbb{Z}_{3} \times \mathbb{Z}_{13}$.}
    \label{B:s4p2}
    \end{figure}

We visualise $\PP$ as a set of $e_4$ copies of the three-dimensional grid $\PP' := \mathbb{Z}_{e_1} \times \mathbb{Z}_{e_2} \times \mathbb{Z}_{e_3}$. We see, from the coordinates listed in Table \ref{B:s4p2-coord}, that the $240$ points of $B_4$ are distributed among the $240$ grids $\PP' \times \{\delta_4\}$ with $1 \leq \delta_4 \leq 240$, such that each $\PP' \times \{\delta_4\}$ contains a unique point of $B_4$. To save on space, we illustrate $B$ by projecting it onto $\PP'$, as in Figure \ref{B:s4p2}, and for clarity we show $\PP'$ as a set of $e_3$ two-dimensional grids $\mathbb{Z}_{e_1} \times \mathbb{Z}_{e_2}$. The points in $B_1 \cup B_2 \cup B_3$ all belong to $\PP' \times \{0\}$; the subsets $B_1$ and $B_2$ in the layer $\delta_3 = 0$ are as labelled, and the square dots in the layers $\delta_3 =1,\dots, 12$ denote the points of $B_3$. The remaining unlabelled circular dots are the points of $B_4$; each of these dots represents two points of $B_4$ that have the same $(\delta_1,\delta_2,\delta_3)$-coordinates but which belong to different grids $\PP' \times \{\delta_4\}$.

The set $B_4$ is divided into subsets which are shown in the three diagrams in the bottom row of Figure \ref{B:s4p2}. The subsets $B_{4,i}$, $1 \leq i \leq 6$, are defined according to the $(\delta_1,\delta_2)$-coordinates (with no restriction on $\delta_3$ and $\delta_4$), and $B_4$ is the disjoint union of these subsets. We also divide the points of $B_4$ with $\delta_2 > 1$ into subsets $B_{4,7}$ and $B_{4,8}$ according to the $(\delta_2,\delta_3)$-coordinates, with no restriction on $\delta_1$ and $\delta_4$. Note that the points in $B_{4,7}$ are those which have the same $(\delta_2,\delta_3)$-coordinates as a point in $B_3$; indeed, in the figure the subset $B_{4,7}$ consists of the points which belong to the same $\delta_3$-layer and to the same row as a point in $B_3$. The points in $B_{4,8}$ are those which do not have this property.

In view of Proposition \ref{lem:B'}, to verify that $\Des{4}{2}$ is a $2$-design we only need to show that condition \eqref{eq:2des} holds for each nonempty subset $J$ of $I = \{1,2,3,4\}$ which does not contain $4$. The array values $x_{\delta_J}$, for each such subset $J$ and each $\delta \in B$, are listed in Table \ref{array:s4p2}, and these can be deduced from Figure \ref{B:s4p2} and from the definition of $B$. Thus each of these parameters can be interpreted as the following numbers in Figure \ref{B:s4p2}:
    \begin{itemize}
    \item $x_{(\delta_1,\delta_2,\delta_3)} =$ number of points that appear in the same $\delta_3$-layer as $\delta$;
    \item $x_{(\delta_2,\delta_3)} =$ number of points that appear in the same row and  same $\delta_3$-layer as $\delta$;
    \item $x_{(\delta_1,\delta_3)} =$ number of points that appear in the same column and same $\delta_3$-layer as $\delta$;
    \item $x_{(\delta_1,\delta_2)} =$ number of points that appear in the same column and the same row, and in any layer, as $\delta$;
    \item $x_{(\delta_3)} =$ number of points that appear in the same $\delta_3$-layer as $\delta$;
    \item $x_{(\delta_2)} =$ number of points that appear in the same row, and in any layer, as $\delta$;
    \item $x_{(\delta_1)} =$ number of points that appear in the same column, and in any layer, as $\delta$.
    \end{itemize}

    \begin{table}[ht]
        \centering
        \begin{tabular}[t]{|crr|}
        \hline
        $\bm{\delta} = \delta|_{\{1,2,3\}}$ & $\# \bm{\delta}$ & $x_{\bm{\delta}}$ \\
        \hline
        $\delta \in B' \times \{0\}$ & $17$ & $1$ \\
        $\delta \in B_4$ & $120$ & $2$ \\
        \hline
        \end{tabular}
        \begin{tabular}[t]{|crr|}
        \hline
        $\delta_3$ & $\# \delta_3$ & $x_{(\delta_3)}$ \\
        \hline
        $0$ & $1$ & $5$ \\
        $\delta_3 \geq 1$ & $12$ & $21$ \\
        \hline
        \end{tabular}
        \begin{tabular}[t]{|crr|}
        \hline
        $\delta_2$ & $\# \delta_2$ & $x_{(\delta_2)}$ \\
        \hline
        $0$ & $1$ & $75$ \\
        $1,2$ & $2$ & $91$ \\
        \hline
        \end{tabular}
        \begin{tabular}[t]{|crr|}
        \hline
        $\delta_1$ & $\# \delta_1$ & $x_{(\delta_1)}$ \\
        \hline
        $0,1,2$ & $3$ & $33$ \\
        $3$ & $1$ & $50$ \\
        $4,5,6$ & $3$ & $36$ \\
        \hline
        \end{tabular}
        \begin{tabular}[t]{|crr|}
        \hline
        $\bm{\delta} = \delta|_{\{1,3\}}$ & $\# \bm{\delta}$ & $x_{\bm{\delta}}$ \\
        \hline
        $\delta \in B_1$ & $3$ & $1$ \\
        $\delta \in B_2$ & $1$ & $2$ \\
        $\delta \in B_3$ & $12$ & $1$ \\
        $\delta \in B_4$ & $60$ & $4$ \\
        \hline
        \end{tabular}
        \quad
        \begin{tabular}[t]{|crr|}
        \hline
        $\bm{\delta} = \delta|_{\{2,3\}}$ & $\# \bm{\delta}$ & $x_{\bm{\delta}}$ \\
        \hline
        $\delta \in B_1$ & $1$ & $3$ \\
        $\delta \in B_2$ & $2$ & $1$ \\
        $\delta \in B_3 \cup B_{4,7}$ & $12$ & $5$ \\
        $\delta \in \bigcup_{i \in \{1,3,5\}} B_{\{4,i\}}$ & $12$ & $6$ \\
        $\delta \in B_{4,8}$ & $12$ & $10$ \\
        \hline
        \end{tabular}
        \quad
        \begin{tabular}[t]{|crr|}
        \hline
        $\bm{\delta} = \delta|_{\{1,2\}}$ & $\# \bm{\delta}$ & $x_{\bm{\delta}}$ \\
        \hline
        $\delta \in B_1 \cup B_{\{4,1\}}$ & $3$ & $9$ \\
        $\delta \in B_{\{4,2\}}$ & $6$ & $12$ \\
        $\delta \in B_{\{4,3\}}$ & $1$ & $24$ \\
        $\delta \in B_2 \cup B_{\{4,4\}}$ & $2$ & $13$ \\
        $\delta \in B_{\{4,5\}}$ & $3$ & $8$ \\
        $\delta \in B_3 \cup B_{\{4,6\}}$ & $6$ & $14$ \\
        \hline
        \end{tabular}
        \caption{Array values for $B$ in Example \ref{ex:s4p2}}
        \label{array:s4p2}
    \end{table}

Arguing as in the proof of Proposition \ref{p:Des(3,p)} in the three-dimensional case, the sum $\sum_{\delta \in B} x_{\delta_J}$ can be obtained from Table \ref{array:s4p2} by taking the sum of the products of the entry in the column ``$\# \bm{\delta}$'' and the square of the corresponding entry in the column ``$x_{\delta_J}$.'' For example, for $J = \{1,2,3\}$,
    \[
    \sum_{(\delta_1,\delta_2,\delta_3) \in \mathbb{Z}_7 \times \mathbb{Z}_3 \times \mathbb{Z}_{13}} x^2_{(\delta_1,\delta_2,\delta_3)}
    = 17 \cdot 1^2 + 120 \cdot 2^2
    = 497
    = 257 + 240
    = k + e_4 - 1.
    \]
So condition \eqref{eq:2des} holds for $J = \{1,2,3\}$. The other conditions can be verified in a similar way, thus proving that $\Des{4}{2}$ is a $2$-design.

\section{Acknowledgements}

The first and third authors (Alavi and Daneshkhah) are grateful to Alice Devillers and Cheryl E. Praeger for supporting their visit to UWA (The University of Western Australia) during February-June 2023. They also thank Bu-Ali Sina University for the support during their sabbatical leave. 

This work forms part of the Australian Research Council Discovery Grant project DP200100080.

\end{document}

%% file: pictures-gridimprim.tex

\def\twogridsquare{
    \multido{\n=-0.75+0.50}{4}{\psline[linecolor=gray](\n,-0.85)(\n,0.85) \psline[linecolor=gray](-0.85,\n)(0.85,\n)}
    \psframe[linecolor=black](-1,-1)(1,1)
    {\color{gray}
        \rput[b](-0.75,1.1){\small $0$} \rput[b](-0.25,1.1){\small $1$} \rput[b](0.25,1.1){\small $2$} \rput[b](0.75,1.1){\small $3$}
        \rput[r](-1.15,0.75){\small $0$} \rput[r](-1.15,0.25){\small $1$} \rput[r](-1.15,-0.25){\small $2$} \rput[r](-1.15,-0.75){\small $3$}}
    }

\def\exampleblocks{
    \rput(-2,0){\twogridsquare \qdisk(-0.75,0.75){3pt} \qdisk(-0.75,0.25){3pt} \qdisk(-0.25,0.25){3pt} \qdisk(-0.25,-0.25){3pt} \qdisk(0.25,-0.25){3pt} \rput[t](0,-1.2){$B$}}
    \rput(2,0){\twogridsquare \qdisk(-0.75,0.25){3pt} \qdisk(-0.75,-0.25){3pt} \qdisk(-0.25,0.25){3pt} \qdisk(-0.25,-0.25){3pt} \qdisk(0.25,0.75){3pt} \rput[t](0,-1.2){$B'$}}
    }


\def\twogridparts{
    \rput(-3,0){
        \multido{\n=-0.75+0.50}{4}{\psline[linecolor=lightgray](-1.6,\n)(1.6,\n)}
        \multido{\n=-1.5+0.5}{7}{\psline[linewidth=1.5pt](\n,-0.85)(\n,0.85)}
        \psframe(-1.75,-1)(1.75,1)
        \rput[b](0,1.25){$\C_{\{1\}}$} \rput[t](0,-1.25){$e_1$ columns}
    }
    \rput(3,0){
        \multido{\n=-1.5+0.5}{7}{\psline[linecolor=lightgray](\n,-0.85)(\n,0.85)}
        \multido{\n=-0.75+0.50}{4}{\psline[linewidth=1.5pt](-1.6,\n)(1.6,\n)}
        \psframe(-1.75,-1)(1.75,1)
        \rput[b](0,1.25){$\C_{\{2\}}$} \rput[t](0,-1.25){$e_2$ rows}
    }
    }


\def\reddots{\multido{\n=0.0+0.5}{2}{\pscircle*[linecolor=red](0,\n){2.5pt}}}

\def\twogridB{
    \multido{\n=-3.0+0.5}{13}{\psline[linecolor=gray](\n,-1.6)(\n,1.6)}
    \multido{\n=-1.5+0.5}{7}{\psline[linecolor=gray](-3.1,\n)(3.1,\n)}
    \psframe(-3.25,-1.75)(3.25,1.75)
    \multido{\n=-3.0+0.5}{4}{\pscircle*[linecolor=amber](\n,1.5){2.5pt}}
    \rput(0,-1.5){\reddots} \rput(-0.5,-0.5){\reddots} \rput(-1,0.5){\reddots}
    \psset{arm=5pt,nodesep=2pt}
    \pnode(-3,1.75){x1} \pnode(-1.5,1.75){x2} \pnode(-1,1.75){x3} \pnode(0,1.75){x4} \pnode(-3.25,-1.75){x5} \pnode(3.25,-1.75){x6}
    \pnode(3.25,-1.5){y1} \pnode(3.25,1){y2} \pnode(-3.25,1.75){y3}
    \ncbar[angle=90,linecolor=amber]{x1}{x2} \naput{\small\color{amber} $p+1$}
    {\psset{linecolor=red} \small\color{red}
        \ncbar[angle=90]{x3}{x4} \naput{$\frac{p^2 - p}{2}$ cols.}
        \ncbar{y1}{y2} \nbput{$p^2 - p$ rows}}
    \rput(-2.25,1){\psframebox*[linecolor=white]{\small\color{amber} $B_1$}}
    \rput(0,0.25){\psframebox*[linecolor=white]{\small\color{red} $B_2$}}
    {\psset{linecolor=gray} \small\color{gray}
        \ncbar[angle=-90]{x6}{x5} \naput{$e_1$ columns}
        \ncbar[angle=180]{x5}{y3} \naput{$e_2$ rows}}
    }


\def\grid{
    \multido{\n=-0.5+0.5}{3}{\pstThreeDPut(0,\n,0){
        \multido{\n=-0.5+0.5}{3}{\pstThreeDLine[linecolor=lightgray](\n,0,-0.5)(\n,0,0.5)}
        }}
    \multido{\n=-0.5+0.5}{3}{\pstThreeDPut(0,0,\n){
        \multido{\n=-0.5+0.5}{3}{
            \pstThreeDLine[linecolor=lightgray](\n,-0.5,0)(\n,0.5,0)
            \pstThreeDLine[linecolor=lightgray](-0.5,\n,0)(0.5,\n,0)
        }}}}
\def\pointset{
    \multido{\n=-0.5+0.5}{3}{\pstThreeDPut(0,0,\n){
        \multido{\n=-0.5+0.5}{3}{\pstThreeDPut(0,\n,0){
            \multido{\n=-0.5+0.5}{3}{\pstThreeDDot(\n,0,0)}
            }}}}}
\def\labels{
    \pstThreeDLine(0.5,0.65,-0.5)(0.5,0.75,-0.5)(-0.5,0.75,-0.5)(-0.5,0.65,-0.5) \pstThreeDPut(0,1.25,-0.5){$e_1$}
    \pstThreeDLine(0.6,-0.5,-0.5)(0.7,-0.5,-0.5)(0.7,0.5,-0.5)(0.6,0.5,-0.5) \pstThreeDPut(0.95,0,-0.5){$e_2$}
    \pstThreeDLine(0.5,-0.7,-0.5)(0.5,-0.8,-0.5)(0.5,-0.8,0.5)(0.5,-0.7,0.5) \pstThreeDPut(0.5,-1.2,0){$e_3$}
    }
\def\labelsII{
    \pstThreeDLine(0.5,0.6,-0.5)(0.5,0.7,-0.5)(-0.5,0.7,-0.5)(-0.5,0.6,-0.5) \pstThreeDPut(0,0.9,-0.5){$e_1$}
    \pstThreeDLine(0.65,-0.5,-0.5)(0.7,-0.5,-0.5)(0.75,0.5,-0.5)(0.65,0.5,-0.5) \pstThreeDPut(1.15,0,-0.5){$e_2$}
    \pstThreeDLine(0.5,-0.6,-0.5)(0.5,-0.7,-0.5)(0.5,-0.7,0.5)(0.5,-0.6,0.5) \pstThreeDPut(0.5,-0.9,0){$e_3$}
    }

\def\threegridlines{
    \psset{Alpha=22,Beta=15,unit=2}
    \rput(-2.5,0){
        \grid
        \multido{\n=-0.5+0.5}{3}{\pstThreeDPut(0,0,\n){
            \multido{\n=-0.5+0.5}{3}{\pstThreeDPut(0,\n,0){
                \pstThreeDLine[linewidth=1.5pt](-0.5,0,0)(0.5,0,0)}
                }}}
        \pointset \labels \rput(0,1){$\C_{\{2,3\}}$}
        }
    {\psset{Alpha=68}
        \grid
        \multido{\n=-0.5+0.5}{3}{\pstThreeDPut(0,0,\n){
            \multido{\n=-0.5+0.5}{3}{\pstThreeDPut(\n,0,0){
                \pstThreeDLine[linewidth=1.5pt](0,-0.5,0)(0,0.5,0)}
                }}}
        \pointset \labelsII \rput(0,1){$\C_{\{1,3\}}$}}
    \rput(2.5,0){
        \grid
        \multido{\n=-0.5+0.5}{3}{\pstThreeDPut(0,\n,0){
            \multido{\n=-0.5+0.5}{3}{\pstThreeDPut(\n,0,0){
                \pstThreeDLine[linewidth=1.5pt](0,0,-0.5)(0,0,0.5)}
                }}}
        \pointset \labels \rput(0,1){$\C_{\{1,2\}}$}}
    }

\def\threegridplanes{
    \psset{Alpha=22,Beta=15,unit=2}
    \rput(-2.5,0){
        \multido{\n=-0.5+0.5}{3}{\pstThreeDPut(\n,0,0){
            \pscustom[fillstyle=solid,fillcolor=gray!20!white]{
                \pstThreeDLine(0,-0.5,-0.5)(0,-0.5,0.5)(0,0.5,0.5)(0,0.5,-0.5)(0,-0.5,-0.5)
                }}}
        \grid \pointset \labels \rput(0,1){$\C_{\{1\}}$}
        \multido{\n=-0.5+0.5}{3}{\pstThreeDPut(\n,0,0){
            \pstThreeDLine(0,-0.5,-0.5)(0,-0.5,0.5)(0,0.5,0.5)(0,0.5,-0.5)(0,-0.5,-0.5)
            }}}
    {\psset{Alpha=68}
        \multido{\n=-0.5+0.5}{3}{\pstThreeDPut(0,\n,0){
            \pscustom[linecolor=gray!20!white,fillstyle=solid,fillcolor=gray!20!white]{
                \pstThreeDLine(-0.5,0,-0.5)(-0.5,0,0.5)(0.5,0,0.5)(0.5,0,-0.5)(-0.5,0,-0.5)
                }}}
        \grid \pointset \labelsII \rput(0,1){$\C_{\{2\}}$}
        \multido{\n=-0.5+0.5}{3}{\pstThreeDPut(0,\n,0){
            \pstThreeDLine(-0.5,0,-0.5)(-0.5,0,0.5)(0.5,0,0.5)(0.5,0,-0.5)(-0.5,0,-0.5)
            }}}
    \rput(2.5,0){
        \multido{\n=-0.5+0.5}{3}{\pstThreeDPut(0,0,\n){
            \pscustom[linecolor=gray!20!white,fillstyle=solid,fillcolor=gray!15!white]{
                \pstThreeDLine(-0.5,-0.5,0)(0.5,-0.5,0)(0.5,0.5,0)(-0.5,0.5,0)(-0.5,-0.5,0)
                }}}
        \grid \pointset \labels \rput(0,1){$\C_{\{3\}}$}
        \multido{\n=-0.5+0.5}{3}{\pstThreeDPut(0,0,\n){
            \pstThreeDLine(-0.5,-0.5,0)(0.5,-0.5,0)(0.5,0.5,0)(-0.5,0.5,0)(-0.5,-0.5,0)
            }}}
    }


\def\twogrid{
    \multido{\n=-1.5+0.5}{7}{\psline[linecolor=gray](\n,-0.65)(\n,0.65)}
    \multido{\n=-0.5+0.5}{3}{\psline[linecolor=gray](-1.65,\n)(1.65,\n)}
    }
\def\twogridcol{\rput[t](0,-1){\footnotesize \color{gray}
    \rput(-1.5,0){$0$} \rput(-1,0){$1$} \rput(-0.5,0){$2$} \rput(0,0){$3$} \rput(0.5,0){$4$} \rput(1,0){$5$} \rput(1.5,0){$6$}
    }}
\def\twogridrow{\rput[r](-1.95,0){\footnotesize \color{gray}
    \rput(0,0.5){$0$} \rput(0,0){$1$} \rput(0,-0.5){$2$}
    }}
\def\threeblock{
    \twogrid \psframe(-1.75,-0.75)(1.75,0.75) \twogridrow \twogridcol
    }

\def\threegridB{
    {\psset{linecolor=amber} \multido{\n=-1.5+0.5}{3}{\qdisk(\n,0.5){2.5pt}}}
    {\psset{linecolor=red} \multido{\n=-0.5+0.5}{2}{\qdisk(0,\n){2.5pt}}}
    {\psset{linecolor=blue} \multido{\n=0.5+0.5}{3}{\rput(\n,0){\multido{\n=-0.5+0.5}{2}{\qdisk(0,\n){2.5pt}}}}}
    \psset{arm=5pt,nodesep=2pt}
    \pnode(-1.5,0.75){x1} \pnode(-0.5,0.75){x2} \pnode(-0.1,-0.55){y1} \pnode(-0.1,0.05){y2}
    \ncbar[angle=90,linecolor=amber]{x1}{x2} \naput{\small\color{amber} $B_1 \times \{(0,0)\}$}
    \ncbar[angle=180,linecolor=red]{y1}{y2} \naput{\psframebox*[linecolor=white]{\small\color{red} $B_2 \times \{0\}$}}
    \rput(1,0.4){\psframebox*[linecolor=white]{\small\color{blue} $B_3$}}
    }


\def\gridI{
    \psline[linecolor=gray](0,-3.5)(0,4)
    \psline[linecolor=gray](-3,-3.5)(-3,4)
    \psline[linecolor=gray](-5.5,-3.5)(-5.5,4)
    \psline[linecolor=gray](2.5,-3.5)(2.5,4)
    \psline[linecolor=gray](3,-3.5)(3,4)
    \psline[linecolor=gray](-6,0)(5.5,0)
    \psline[linecolor=gray](-6,3.5)(5.5,3.5)
    \psframe(-6,-3.5)(5.5,4)
    {\psset{linecolor=amber,fillcolor=yellow!50!white,fillstyle=solid}
    \psframe(-6,4)(-5.5,3.5) \psframe(-5.5,4)(-3,3.5)
    }
    {\psset{linecolor=blue,fillcolor=lightblue,fillstyle=solid}
    \psframe(-3,3.5)(0,0) 
    \psframe(2.5,0)(3,3.5) 
    \psframe(3,0)(5.5,-3.5) 
    \rput(0,-3.5){\psframe(-3,3.5)(-2.7,2.8) \psframe(-2.7,2.8)(-2.4,2.1) \psframe(-2.4,2.1)(-2.1,1.4) \psframe(-2.1,1.4)(-1.8,0.7) \psframe(-1.8,0.7)(-1.5,0)} 
    \rput(-7.5,-1.75){\psframe(0,0.35)(0.3,-0.35) \rput(0.6,0){$=$} \psline[linecolor=lightgray](1,-0.5)(1,0.5) \pscircle*(1,0.25){2.5pt} \pscircle*(1,-0.25){2.5pt}}
    }
    \psline[linecolor=gray,linestyle=dashed](-1.5,-3.5)(-1.5,3.5)
    {\psset{linecolor=red,fillcolor=pink,fillstyle=solid}
    \psframe(-3,3.5)(-2.7,2.8) \psframe(-2.7,2.8)(-2.4,2.1) \psframe(-2.4,2.1)(-2.1,1.4) \psframe(-2.1,1.4)(-1.8,0.7) \psframe(-1.8,0.7)(-1.5,0) \psframe(-1.5,0)(-1.2,-0.7) \psframe(-1.2,-0.7)(-0.9,-1.4) \psframe(-0.9,-1.4)(-0.6,-2.1) \psframe(-0.6,-2.1)(-0.3,-2.8) \psframe(-0.3,-2.8)(0,-3.5)
    \rput(-7.5,3){\psframe(0,0.35)(0.3,-0.35) \rput(0.6,0){$=$} \psline[linecolor=lightgray](1,-0.5)(1,0.5) \pscircle*(1,0.25){2.5pt} \pscircle*(1,-0.25){2.5pt}}
    }
    \pnode(5.5,4){y1} \pnode(5.5,3.5){y2} 
    \ncbar[nodesep=2pt,arm=3pt,linewidth=0.25pt]{y1}{y2} \naput{\small $1$ row}
    {\psset{nodesep=2pt,arm=3pt,linewidth=0.25pt,angle=-90}
    \pnode(-6,-3.5){x1} \pnode(-5.5,-3.5){x2} \pnode(-3,-3.5){x3} \pnode(-1.5,-3.5){x4} \pnode(0,-3.5){x5} \pnode(2.5,-3.5){x6} \pnode(3,-3.5){x7} \pnode(5.5,-3.5){x8}
    \ncbar{x1}{x2} \nbput{\small \parbox{0.6cm}{\centering $1$ \\ col}}
    \ncbar{x2}{x3} \nbput{\small $p$ cols}
    \ncbar{x3}{x4} \nbput{\small \parbox{2cm}{\centering $\frac{p^2 - p}{4}$ \\ cols}}
    \ncbar{x4}{x5} \nbput{\small \parbox{2cm}{\centering $\frac{p^2 - p}{4}$ \\ cols}}
    \ncbar{x5}{x6} \nbput{\small $p$ cols}
    \ncbar{x6}{x7} \nbput{\small \parbox{0.6cm}{\centering $1$ \\ col}}
    \ncbar{x7}{x8} \nbput{\small \parbox{2cm}{\centering $\frac{p^2 - p}{2} - 1$ \\ cols}}
    }
    \rput[b](-5.75,4.1){\small\color{amber} $B_{1,1}$} \rput[b](-4.25,4.1){\small\color{amber} $B_{1,2}$}
    \rput(-0.5,-0.5){\small\color{red} $B_2$}
    \rput(-2.55,-2.45){\small\color{blue} $B_{3,2}$} \rput(-2.5,0.5){\small\color{blue} $B_{3,2}$} \rput(-1,1.8){\small\color{blue} $B_{3,2}$}
    \rput[l](3.2,1.75){\small\color{blue} $B_{3,4}$} \rput(4.25,-1.75){\small\color{blue} $B_{3,4}$}
}

\def\gridII{
    \psline[linecolor=gray](0,-4.2)(0,4)
    \psline[linecolor=gray](-3.3,-4.2)(-3.3,4)
    \psline[linecolor=gray](-5.8,-4.2)(-5.8,4)
    \psline[linecolor=gray](2.5,-4.2)(2.5,4)
    \psline[linecolor=gray](3,-4.2)(3,4)
    \psline[linecolor=gray](-6.3,-0.35)(5.5,-0.35)
    \psline[linecolor=gray](-6.3,3.5)(5.5,3.5)
    \psframe(-6.3,-4.2)(5.5,4)
    {\psset{linecolor=amber,fillcolor=yellow!50!white,fillstyle=solid}
    \psframe(-6.3,4)(-5.8,3.5) \psframe(-5.8,4)(-3.3,3.5) 
    }
    {\psset{linecolor=blue,fillcolor=lightblue,fillstyle=solid}
    \psframe(-3.3,3.5)(0,-0.35) 
    \psframe(2.5,3.5)(3,-0.35) 
    \psframe(3,-0.35)(5.5,-4.2) 
    \rput(-0.3,-3.85){\psframe(-3,3.5)(-2.7,2.8) \psframe(-2.7,2.8)(-2.4,2.1) \psframe(-2.4,2.1)(-2.1,1.4) \psframe(-2.1,1.4)(-1.8,0.7) \psframe(-1.8,0.7)(-1.5,0)} \psframe(-1.8,-3.85)(-1.5,-4.2)
    \rput(-7.8,-1.75){\psframe(0,0.35)(0.3,-0.35) \rput(0.6,0){$=$} \psline[linecolor=lightgray](1,-0.5)(1,0.5) \pscircle*(1,0.25){2.5pt} \pscircle*(1,-0.25){2.5pt}}
    }
    \psline[linecolor=gray,linestyle=dashed](-1.5,-4.2)(-1.5,3.5)
    {\psset{linecolor=red,fillcolor=pink,fillstyle=solid}
    \rput(-0.3,0){
    \psframe(-3,3.5)(-2.7,2.8) \psframe(-2.7,2.8)(-2.4,2.1) \psframe(-2.4,2.1)(-2.1,1.4) \psframe(-2.1,1.4)(-1.8,0.7) \psframe(-1.8,0.7)(-1.5,0) \psframe(-1.5,0)(-1.2,-0.35) \psframe(-1.5,-0.35)(-1.2,-0.7) \psframe(-1.2,-0.7)(-0.9,-1.4) \psframe(-0.9,-1.4)(-0.6,-2.1) \psframe(-0.6,-2.1)(-0.3,-2.8) \psframe(-0.3,-2.8)(0,-3.5)} \psframe(-0.3,-3.5)(0,-4.2)
    \rput(-7.8,3){\psframe(0,0.35)(0.3,-0.35) \rput(0.6,0){$=$} \psline[linecolor=lightgray](1,-0.5)(1,0.5) \pscircle*(1,0.25){2.5pt} \pscircle*(1,-0.25){2.5pt}}
    }
    \pnode(5.5,4){y1} \pnode(5.5,3.5){y2}
    \ncbar[nodesep=2pt,arm=3pt,linewidth=0.25pt]{y1}{y2} \naput{\small $1$ row}
    {\psset{nodesep=2pt,arm=3pt,linewidth=0.25pt,angle=-90}
    \pnode(-6.3,-4.2){x1} \pnode(-5.8,-4.2){x2} \pnode(-3.3,-4.2){x3} \pnode(-1.5,-4.2){x4} \pnode(0,-4.2){x5} \pnode(2.5,-4.2){x6} \pnode(3,-4.2){x7} \pnode(5.5,-4.2){x8}
    \ncbar{x1}{x2} \nbput{\small \parbox{0.6cm}{\centering $1$ \\ col}}
    \ncbar{x2}{x3} \nbput{\small $p$ cols}
    \ncbar{x3}{x4} \nbput{\small \parbox{2cm}{\centering $\frac{p^2 - p + 2}{4}$ \\ cols}}
    \ncbar{x4}{x5} \nbput{\small \parbox{2cm}{\centering $\frac{p^2 - p - 2}{4}$ \\ cols}}
    \ncbar{x5}{x6} \nbput{\small $p$ cols}
    \ncbar{x6}{x7} \nbput{\small \parbox{0.6cm}{\centering $1$ \\ col}}
    \ncbar{x7}{x8} \nbput{\small \parbox{2cm}{\centering $\frac{p^2 - p}{2} - 1$ \\ cols}}
    }
    \rput[b](-6.05,4.1){\small\color{amber} $B_{1,1}$} \rput[b](-4.55,4.1){\small\color{amber} $B_{1,2}$}
    \rput(-2.8,-3){\small\color{blue} $B_{3,2}$} \rput(-2.7,0.25){\small\color{blue} $B_{3,2}$} \rput(-1,1.58){\small\color{blue} $B_{3,2}$}
    \rput[l](3.2,1.75){\small\color{blue} $B_{3,4}$} \rput(4.25,-2.4){\small\color{blue} $B_{3,4}$}
}


\def\threegridBgenImodIV{
    \psline[linecolor=gray](-6,-1.4)(5.5,-1.4)
    \gridI
    {\psset{linecolor=blue,fillcolor=lightblue,fillstyle=solid}
    \psframe(-5.5,3.5)(-3,-1.4) 
    \psframe(0,-3.5)(2.5,-1.4) 
    }
    {\psset{linecolor=gray,linestyle=dashed}
    \psline(-5.5,0)(-3,0) \psline(3,-1.4)(5.5,-1.4)
    }
    {\psset{nodesep=2pt,arm=3pt,linewidth=0.25pt}
    \pnode(5.5,0){y3} \pnode(5.5,-1.4){y4} \pnode(5.5,-3.5){y5}
    \ncbar{y2}{y3} \naput{\small $\frac{p^2 - p}{2}$ rows}
    \ncbar{y3}{y4} \naput{\small $\frac{p-1}{2}$ rows}
    \ncbar{y4}{y5} \naput{\small $\frac{p^2 - 2p + 1}{2}$ rows}
    }
    \rput(-4.25,1.05){\small\color{blue} $B_{3,1}$}
    \rput(1.25,-2.45){\small\color{blue} $B_{3,3}$}
    }


\def\threegridBgenIIImodIV{
    \psline[linecolor=gray](-6.3,-2.1)(5.5,-2.1)
    \gridII
    {\psset{linecolor=blue,fillcolor=lightblue,fillstyle=solid}
    \psframe(-5.8,3.5)(-3.3,-2.1) 
    \psframe(0,-2.1)(2.5,-4.2) 
    }
    {\psset{linecolor=gray,linestyle=dashed}
    \psline(-5.8,-0.35)(-3.3,-0.35) \psline(3,-2.1)(5.5,-2.1)
    }
    {\psset{nodesep=2pt,arm=3pt,linewidth=0.25pt}
    \pnode(5.5,-0.35){y3} \pnode(5.5,-2.1){y4} \pnode(5.5,-4.2){y5}
    \ncbar{y2}{y3} \naput{\small $\frac{p^2 - p}{2}$ rows}
    \ncbar{y3}{y4} \naput{\small $\frac{p-1}{2}$ rows}
    \ncbar{y4}{y5} \naput{\small $\frac{p^2 - 2p + 1}{2}$ rows}
    }
    \rput(-4.55,0.7){\small\color{blue} $B_{3,1}$}
    \rput(1.25,-2.8){\small\color{blue} $B_{3,3}$}
    }


\def\threegridBgenOmodIV{
    \psline[linecolor=gray](-6,1.4)(5.5,1.4)
    \gridI
    {\psset{linecolor=blue,fillcolor=lightblue,fillstyle=solid}
    \psframe(-5.5,3.5)(-3,1.4) 
    \psframe(0,-3.5)(2.5,1.4) 
    }
    {\psset{linecolor=gray,linestyle=dashed}
    \psline(-3,1.4)(0,1.4) \psline(2.5,1.4)(3,1.4)
    }
    {\psset{nodesep=2pt,arm=3pt,linewidth=0.25pt}
    \pnode(5.5,1.4){y3} \pnode(5.5,0){y4} \pnode(5.5,-3.5){y5}
    \ncbar{y2}{y3} \naput{\small $\frac{p^2 - 2p}{2}$ rows}
    \ncbar{y3}{y4} \naput{\small $\frac{p}{2}$ rows}
    \ncbar{y4}{y5} \naput{\small $\frac{p^2 - p}{2}$ rows}
    }
    \rput(-4.25,2.45){\small\color{blue} $B_{3,1}$}
    \rput(1.25,-1.05){\small\color{blue} $B_{3,3}$}
    }

    
\def\threegridBgenIImodIV{
    \psline[linecolor=gray](-6.3,1.4)(5.5,1.4)
    \gridII
    {\psset{linecolor=blue,fillcolor=lightblue,fillstyle=solid}
    \psframe(-5.8,3.5)(-3.3,1.4)
    \psframe(0,1.4)(2.5,-4.2)
    }
    {\psset{linecolor=gray,linestyle=dashed}
    \psline(-3.3,1.4)(0,1.4) \psline(2.5,1.4)(3,1.4) \psline(0,-0.35)(2.5,-0.35)
    }
    {\psset{nodesep=2pt,arm=3pt,linewidth=0.25pt}
    \pnode(5.5,1.4){y3} \pnode(5.5,-0.35){y4} \pnode(5.5,-4.2){y5}
    \ncbar{y2}{y3} \naput{\small $\frac{p^2 - 2p}{2}$ rows}
    \ncbar{y3}{y4} \naput{\small $\frac{p}{2}$ rows}
    \ncbar{y4}{y5} \naput{\small $\frac{p^2 - p}{2}$ rows}
    }
    \rput(-4.55,2.45){\small\color{blue} $B_{3,1}$}
    \rput(1.25,-1.57){\small\color{blue} $B_{3,3}$}
    }


\def\fourgridBI{
    \threeblock \rput[b](0,0.95){\footnotesize \color{gray} $\delta_3 = 0$}
        {\psset{linecolor=amber}\qdisk(-1.5,0.5){2.5pt} \qdisk(-1,0.5){2.5pt} \qdisk(-0.5,0.5){2.5pt}}
        {\psset{linecolor=red}\qdisk(0,0){2.5pt} \qdisk(0,-0.5){2.5pt}}
        \psset{arm=5pt,nodesep=5pt}
        \pnode(-1.55,0.5){x1} \pnode(-0.45,0.5){x2} \pnode(0,-0.55){y1} \pnode(0,0.05){y2}
        \ncbar[angle=-90,linecolor=amber]{x1}{x2} \nbput{\psframebox*[linecolor=white]{\small\color{amber} $B_1$}}
        \ncbar[linecolor=red]{y1}{y2} \nbput{\psframebox*[linecolor=white]{\small\color{red} $B_2$}}
    \end{pspicture}
    
    \begin{pspicture}(-6,-1.35)(6,1.35)
    \rput(-6.5,0){\threeblock \rput[b](0,0.95){\footnotesize \color{gray} $\delta_3 = 1$} \psdot[linecolor=blue,dotstyle=square*,dotsize=6pt](0.5,0)
    {\psset{linecolor=malachite,dotsize=5pt}
        \psdots*(-1,0.5)(0,0.5)(1,0.5)(-1,-0.5)(0,-0.5)(1,-0.5)
        \psdots*(-0.5,0)(1.5,0)(-0.5,-0.5)(1.5,-0.5)
        }
    }
    \rput(-2.15,0){\threeblock \rput[b](0,0.95){\footnotesize \color{gray} $\delta_3 = 2$} \psdot[linecolor=blue,dotstyle=square*,dotsize=6pt](0.5,0)
    {\psset{linecolor=malachite,dotsize=5pt}
        \psdots*(-1,0.5)(0,0.5)(1,0.5)(-1,-0.5)(0,-0.5)(1,-0.5)
        \psdots*(-0.5,0)(1.5,0)(-0.5,-0.5)(1.5,-0.5)
        }
    }
    \rput(2.15,0){\threeblock \rput[b](0,0.95){\footnotesize \color{gray} $\delta_3 = 3$} \psdot[linecolor=blue,dotstyle=square*,dotsize=6pt](0.5,-0.5)
    {\psset{linecolor=malachite,dotsize=5pt}
        \psdots*(-1,0.5)(0,0.5)(1,0.5)(-1,0)(0,0)(1,0)
        \psdots*(-0.5,0)(1.5,0)(-0.5,-0.5)(1.5,-0.5)
        }
    }
    \rput(6.5,0){\threeblock \rput[b](0,0.95){\footnotesize \color{gray} $\delta_3 = 4$} \psdot[linecolor=blue,dotstyle=square*,dotsize=6pt](0.5,-0.5)
    {\psset{linecolor=malachite,dotsize=5pt}
        \psdots*(-1,0.5)(0,0.5)(1,0.5)(-1,0)(0,0)(1,0)
        \psdots*(-0.5,0)(1.5,0)(-0.5,-0.5)(1.5,-0.5)
        }
    }
    }

\def\fourgridBII{
    \rput(-6.5,0){\threeblock \rput[b](0,0.95){\footnotesize \color{gray} $\delta_3 = 5$} \psdot[linecolor=blue,dotstyle=square*,dotsize=6pt](1,0)
    {\psset{linecolor=malachite,dotsize=5pt}
        \psdots(-0.5,0.5)(0,0.5)(1.5,0.5)(-0.5,-0.5)(0,-0.5)(1.5,-0.5)
        \psdots(-1.5,0)(0.5,0)(-1.5,-0.5)(0.5,-0.5)
        }
    }
    \rput(-2.15,0){\threeblock \rput[b](0,0.95){\footnotesize \color{gray} $\delta_3 = 6$} \psdot[linecolor=blue,dotstyle=square*,dotsize=6pt](1,0)
    {\psset{linecolor=malachite,dotsize=5pt}
        \psdots(-0.5,0.5)(0,0.5)(1.5,0.5)(-0.5,-0.5)(0,-0.5)(1.5,-0.5)
        \psdots(-1.5,0)(0.5,0)(-1.5,-0.5)(0.5,-0.5)
        }
    }
    \rput(2.15,0){\threeblock \rput[b](0,0.95){\footnotesize \color{gray} $\delta_3 = 7$} \psdot[linecolor=blue,dotstyle=square*,dotsize=6pt](1,-0.5)
    {\psset{linecolor=malachite,dotsize=5pt}
        \psdots(-0.5,0.5)(0,0.5)(1.5,0.5)(-0.5,0)(0,0)(1.5,0)
        \psdots(-1.5,0)(0.5,0)(-1.5,-0.5)(0.5,-0.5)
        }
    }
    \rput(6.5,0){\threeblock \rput[b](0,0.95){\footnotesize \color{gray} $\delta_3 = 8$} \psdot[linecolor=blue,dotstyle=square*,dotsize=6pt](1,-0.5)
    {\psset{linecolor=malachite,dotsize=5pt}
        \psdots(-0.5,0.5)(0,0.5)(1.5,0.5)(-0.5,0)(0,0)(1.5,0)
        \psdots(-1.5,0)(0.5,0)(-1.5,-0.5)(0.5,-0.5)
        }
    }
    }

\def\fourgridBIII{
    \rput(-6.5,0){\threeblock \rput[b](0,0.95){\footnotesize \color{gray} $\delta_3 = 9$} \psdot[linecolor=blue,dotstyle=square*,dotsize=6pt](1.5,0)
    {\psset{linecolor=malachite,dotsize=5pt}
        \psdots(-1.5,0.5)(0,0.5)(0.5,0.5)(-1.5,-0.5)(0,-0.5)(0.5,-0.5)
        \psdots(-1,0)(1,0)(-1,-0.5)(1,-0.5)
        }
    }
    \rput(-2.15,0){\threeblock \rput[b](0,0.95){\footnotesize \color{gray} $\delta_3 = 10$} \psdot[linecolor=blue,dotstyle=square*,dotsize=6pt](1.5,0)
    {\psset{linecolor=malachite,dotsize=5pt}
        \psdots(-1.5,0.5)(0,0.5)(0.5,0.5)(-1.5,-0.5)(0,-0.5)(0.5,-0.5)
        \psdots(-1,0)(1,0)(-1,-0.5)(1,-0.5)
        }
    }
    \rput(2.15,0){\threeblock \rput[b](0,0.95){\footnotesize \color{gray} $\delta_3 = 11$} \psdot[linecolor=blue,dotstyle=square*,dotsize=6pt](1.5,-0.5)
    {\psset{linecolor=malachite,dotsize=5pt}
        \psdots(-1.5,0.5)(0,0.5)(0.5,0.5)(-1.5,0)(0,0)(0.5,0)
        \psdots(-1,0)(1,0)(-1,-0.5)(1,-0.5)
        }
    }
    \rput(6.5,0){\threeblock \rput[b](0,0.95){\footnotesize \color{gray} $\delta_3 = 12$} \psdot[linecolor=blue,dotstyle=square*,dotsize=6pt](1.5,-0.5)
    {\psset{linecolor=malachite,dotsize=5pt}
        \psdots(-1.5,0.5)(0,0.5)(0.5,0.5)(-1.5,0)(0,0)(0.5,0)
        \psdots(-1,0)(1,0)(-1,-0.5)(1,-0.5)
        }
    }
    }

\def\fourgridBIV{
    \rput(-4.5,0){
        \psline(-1.75,0.25)(1.75,0.25) \psline(-0.25,-0.75)(-0.25,0.75) \psline(0.25,-0.75)(0.25,0.75)
        \psframe(-1.75,-0.75)(1.75,0.75) \twogridrow \twogridcol
        {\color{darkgreen}
        \rput(-1,0.5){\footnotesize $B_{4,1}$} \rput[b](0,0.8){\rnode{p}{\small $B_{4,3}$}} \pnode(0,0.5){q} \ncline[linecolor=darkgreen]{p}{q} \rput(1,0.5){\footnotesize $B_{4,5}$}
        \rput(-1,-0.25){\small $B_{4,2}$} \rput(1,-0.25){\small $B_{4,6}$} \rput{90}(0,-0.25){\footnotesize $B_{4,4}$}
        }}
    \psline(-1.75,0.25)(1.75,0.25) \psline(-1.75,-0.25)(1.75,-0.25)
    \psframe(-1.75,-0.75)(1.75,0.75) \twogridrow \twogridcol
    {\color{darkgreen} \footnotesize \rput(0,0){$B_{4,7}$} \rput(0,-0.5){$B_{4,8}$}}
    \rput[b](0,0.9){\footnotesize \color{gray}$\delta_3 \in \{1,2,5,6,9,10\}$}
    \rput(4.5,0){
        \psline(-1.75,0.25)(1.75,0.25) \psline(-1.75,-0.25)(1.75,-0.25)
        \psframe(-1.75,-0.75)(1.75,0.75) \twogridrow \twogridcol
        {\color{darkgreen} \footnotesize \rput(0,0){$B_{4,8}$} \rput(0,-0.5){$B_{4,7}$}}
        \rput[b](0,0.9){\footnotesize \color{gray}$\delta_3 \in \{3,4,7,8,11,12\}$}
        }
    }